\documentstyle[bal]{article}

\input amssym1.def

\mathchardef\smallsetminus="2\msbfam@72
\mathchardef\varnothing="0\msbfam@3F
\renewcommand{\setminus}{\smallsetminus}
\renewcommand{\emptyset}{\varnothing}

\def\T{\mathop{\cal T}{}}
\def\AA{\mathop{\cal A}{}}

\def\A{\mathop{\bf A}\nolimits}
\def\x{\mathop{\bf x}\nolimits}
\def\y{\mathop{\bf y}\nolimits}
\def\z{\mathop{\bf z}\nolimits}
\def\e{\mathop{\bf e}\nolimits}
\def\v{\mathop{\bf v}\nolimits}
\def\B{\mathop{\bf A'}\nolimits}

\def\hmi{\hspace{-.04em}}
\def\hmii{\hspace{-.1em}}
\def\dasp{d\hmi'\hmi A\hmi s\hmi p\hmi r\hmi e\hmi m\hmi o\hmi n\hmi t,}
\def\cdcc{,\hmii...,\hmii}
\def\ms{\\ $\mathstrut$}
\def\vpx{\mathop{v\_p(x)}\nolimits}
\def\vpi{\mathop{v\_p(i)}\nolimits}
\def\vpj{\mathop{v\_p(j)}\nolimits}
\def\U{\mathop{\cal U}{}}
\def\equi{\mathop{\Leftrightarrow}}
\def\Xtil{\mathop{\widetilde X}{}}
\def\mto{$\to$}

\def\_#1{\mathop{\hspace{-2pt}^{}_{#1}}}
\def\suml  {\mathop{\sum}   \limits}
\def\cupl  {\mathop{\bigcup}\limits}

\def\cdc{,\ldots,}
\def\on{1,\ldots,n}
\def\om{1,\ldots,m}
\def\ok{1,\ldots,k}
\def\si{s\_i}
\def\sj{s\_j}

\def\a{\mathop{\alpha}\nolimits}
\def\ve{\mathop{\varepsilon}\nolimits}
\def\w{\mathop{\varphi}\nolimits}
\def\vp{\mathop{\varphi}\nolimits}
\def\l{\ell}

\def\R{\mathop{\Bbb R}\nolimits}
\def\N{\mathop{\Bbb N}{}}
\def\dej{d\_0\e^j}

\newtheorem{axiom}[theorem]{Axiom}

\title{%
Characterizations of scoring methods for preference aggregation%
\thanks{This work was supported by Russian Foundation for Basic
Research Grant No.~96-01-01010. Partial research support from the
European Community under Grant ACE-91-R02 is also gratefully
acknowledged. The authors thank Anna Khmelnitskaya for some
bibliographic suggestions
}
}
\author{%
Pavel Yu.\ Chebotarev\institute{1} and
Elena Shamis\institute{1}%
\from
\institute{1}%
Institute of Control Sciences, 65 Profsoyuznaya,\\
Moscow 117997, Russia\\
\email {pchv@rambler.ru}
}

\date{May 31, 1997}

\begin{document}
\maketitle

\begin{abstract}

The paper surveys more than forty characterizations of scoring methods
for preference aggregation and contains one new result. A general
scoring operator is {\it self-consistent\/} if alternative $i$ is
assigned a greater score than $j$ whenever $i$ gets no worse (better)
results of comparisons and its `opponents' are assigned respectively
greater (no smaller) scores than those of $j$. We prove that
self-consistency is satisfied if and only if the application of a
scoring operator reduces to the solution of a homogeneous system of
algebraic equations with a monotone function on the left-hand side.

\amsmos
90A07, 90A10, 90A28, 90-02.


\end{abstract}

\section{Introduction}

Scoring procedures transform profiles of individual preferences over a
set of alternatives into scores of the alternatives. The scores can be
used in themselves or serve as the basis for ranking or choice. For the
present, only a few scoring procedures are endowed with their axiomatic
characterizations. At the same time, a large number of ingenious
procedures are advocated and used in such disciplines as management
science, operations research, psychometrics, applied statistics,
processing of sport tournaments, graph theory, etc. Very few social
choice papers deal with them. The aim of this paper is to take one
circumspect step toward an axiomatic framework for comparing the merits
of these elaborate procedures. As a result, we would like to isolate a
family of scoring procedures that comprises a majority of `reasonable'
procedures (so that the further axioms could be imposed on this
family). Two main approaches are applicable. The first one is to
express the desired properties axiomatically, the second is to gather
the existing procedures and specify their common algebraic form. We use
both, and their results are concordant: A scoring procedure satisfies
the axiom of {\it self-consistency} if and only if it has a {\it
monotone implicit form}.

To circumscribe the variety of axioms that have already been used in
the literature (and which can be adapted for the family of procedures
isolated here), we survey a number of characterizations of scoring
methods.

The paper is organized as follows. After the main notation
(Section~\ref{Nota}), we give a review of some papers that characterize
scoring methods for preference aggregation (Section~\ref{Revi}),
introduce the notions of self-consistency (Section~\ref{Self}) and
monotone implicit form (Section~\ref{MIF}), and prove our theorem
(Section~\ref{Theo}).

\section{Main notation}
\label{Nota}

Let $\A=(A^{(1)}\cdc A^{(m)})$ be a profile of individual preferences
over the set of alternatives $X=\{\on\}$. Here $m$ is the number of
individuals and $A^{(p)},\: p=\om,$ represents the preferences of the
$p$th individual. In classical papers, $A^{(p)}$ are linear orders.
Here we would like to involve several other settings as well, where
equivalencies are allowed or transitivity is not assumed or even
degrees of preference are reported.  Namely, $A^{(p)}=(a_{ij}^p)$ is a
completely defined matrix of paired comparisons, where for $i\ne j$,
$a_{ij}^p\in[0,1]$, $a_{ij}^p+a_{ji}^p=1$ and $a_{ii}^p=0$ (although
$a_{ii}^p=0$ are not used in the paper).  When only strict preferences
are allowed, $a_{ij}^p\in\{0,1\}$, otherwise $a_{ij}^p=a_{ji}^p=1/2$
means equivalence; in some settings $a_{ij}^p\in[0,1]$ is interpreted
as the part of a unit preference or probability that is assigned to $i$
in comparison with $j$ (in turn $a_{ji}^p=1-a_{ij}^p$ is assigned to
$j$).  Transitivity may or may not be assumed.

The constraint $a_{ij}^p+a_{ji}^p=1$ may seem to be unduly restrictive,
but we believe it is not. In this setting, $a_{ij}^p$ and $a_{ji}^p$
are not independent comparison outcomes, but rather two complementary
evaluations of the same outcome. That is why their relation merely
characterizes the method of evaluation. Monotone transformations take
them to other popular types of data, e.g., exponential transformations
$\psi$ give $\psi(a_{ij}^p)\cdot\psi(a_{ji}^p)={\rm const}$, which is
typical of the analytical hierarchy framework. Yet, as was noticed by
an anonymous referee, the theorem in Section~\ref{Theo} is adaptable
for the case with mutually independent $a_{ij}^p$ and $a_{ji}^p$ too.

Letting some restrictions on $A^{(p)}$ (or on $\A$ as a whole) be
imposed, suppose $\AA$ is the set of all admissible preference profiles
with given $n$ (the number of alternatives) and $m$ (the number of
individuals).  In Sections~\ref{Self}--\ref{Theo}, $n$ and $m$ are
fixed.  Let $M=\{\om\}$, $\;X\_i=\{1\cdc i-1,\:i+\on\}$, $\;i=\on$.

A {\it scoring operator\/} or {\it scoring procedure\/} is a function
$\w:\:\AA\to \R^n,$ where $\w(\A)=(s\_1\cdc s\_n)$, $\si$ being
a {\it score\/} (or weight, etc.) assigned to alternative $i$. `Scores'
are often associated with total points, however, a number of papers
consider scores as general weights as in this definition, see, e.g.,
\cite{Lasl,Dani}. Let {\em scoring methods\/} be the generic term for
scoring operators and procedures of ranking/choice based on them.

The review in the following section involves a number of results that
deal with social utility (welfare) functions. There is no essential
difference between them and the scores $s\_i$, however, we do not find
it natural to give them a common designation. So, we hold $s\_i$ for
scores derived from paired comparisons and $v(i)$ (or $v(x),$ $\;x\in
X$) for social utilities. Individual utilities will be denoted by
$v\_1(x)\cdc v\_m(x)$, $\;x\in X$.

\section{Axiomatic characterizations of scoring methods}
\label{Revi}

Some papers that characterize scoring methods are presented in
Table~\ref{TCha} (Parts~1 to~6). They are ordered chronologically (up
to years), and every column can be used to classify the results.
Characterizations of scoring procedures that have the form of product
are not included, since their theory is entirely parallel to that of
additive procedures (see, e.g., \cite{GooMar,Fish87,Moul88,BarCooGol}).
We use six standard abbreviations: N for neutrality, A for anonymity, M
for monotonicity, IIA and NIIA for Arrow's and Nash's independences of
irrelevant alternatives, and P for the strong Pareto principle (`weak P'
designates the weak Pareto principle).

\subsection{Resulting preference structures}

Let us start with the third column. The scores of alternatives can be
used in themselves, for ranking or choice. So the results are
classified as ones characterizing:

\begin{itemize}
\item Ranking procedures based on scores. They may provide:
     \begin{itemize}
     \item Weak order\
     \cite{Flem52,Miln54,Fish69,Mork71,Gard73,Smith73,FinFin74,dAspGev77,%
     Mask78,DesGev78,YouLev,Rubi,Robe80,NitRub,Henr,Young86,Bouy92a,Marc96a}.
     \item Weak order over all admissible or feasible vectors of
           individual utilities/scores or on an extension of $X$
           \cite{GooMar,BlaGir54}.
     \item Partial order \cite{FinFin74,BouPer}.
     \end{itemize}
\item Choice procedures based on scores (alternatives with the highest
     score are chosen). Table~\ref{TCha} contains procedures that
     provide:
     \begin{itemize}
     \item Choice set
     \cite{Young74,Young75,HanSah,Fish78,Rich78,FarNit79,Mork82,Young86,%
     SaariSCW,Myer}.
     \item Several variants of choice set with $k$ members \cite{Debo}.
           (It is interesting to compare Debord's result with another
           approach to the committee selection advocated and
           characterized in \cite{ChaCou} and also based on the Borda
           scores.)
     \item Choice from a set of admissible or feasible vectors of
           individual utilities/scores or from an extension of $X$
           \cite{Myer81,Moul88}.
     \item Choice function that attaches a nonempty choice set to every
     subset of $X$ \cite{Fish73,Henr,BaiXu91,Bouy92b,Dugg}.
     \end{itemize}
\item Procedures resulting in scores, measures of social welfare,
      utility, etc., possibly on an extended set of alternatives
      \cite{Hars55,Keen76,Barb,Cheb89a,BarGol,Cheb94,Marc96b}.
\end{itemize}

\hyphenpenalty=10000
\begin{Table}
{Axiomatic characterizations of scoring methods (Part 1)}
\begin{tabular}[t]{p{1.8cm}|p{.1cm}p{2.19cm}p{.2cm}p{1.64cm}p{.2cm}%
p{3.19cm}p{.2cm}p{4.08cm}}
\begin{tabular}{l}
Paper
\end{tabular}
&&\begin{tabular}{l}
Input\\preferences
\end{tabular}
&&\begin{tabular}{l}
Resulting\\(social)\\structure
\end{tabular}
&&\begin{tabular}{l}
Most important\\axioms
\end{tabular}
&&\begin{tabular}{l}
Result of the method
\end{tabular}
\\[2mm]
\hline
\\[-.8mm]
Fleming \cite{Flem52}
&&$m$ weak orders
&& weak order
&&\raggedright
P, independence of locally unaffected individuals
&&{\raggedright
ranking by the sum of arbitrary utilities that represent individual
orders

}
\\[-.8mm]
\raggedright
Goodman, Markowitz \cite{GooMar}
&&\raggedright
A set $\Xtil\subset\R^m$ of admissible utility vectors
$\v=(v\_1\cdcc v\_m)$
(including the discrete case)
&&\raggedright
weak order on $\Xtil$
&&\raggedright
a. P\\
b. independence of the individual zeros, A + `a.'
&&{\raggedright
ranking by:\\
a. $W(v\_1\cdcc v\_m)$ with arbitrary increasing $W:
\R^m\to \R$\\
b. $\sum_{p=1}^mv\_p$
\ms

}
\\[-.8mm]
\raggedright
Blackwell, Girshick \cite{BlaGir54}
&&\raggedright
A set $\Xtil\subset\R^m$ of feasible utility vectors
$\v=(v\_1\cdcc v\_m)$
&&\raggedright
weak order on $\Xtil$
&&\raggedright
strengthened IIA,
nonstrict Pareto preference,
independence of the individual zeros
&&{\raggedright
ranking by
$\sum_{p=1}^m\a\_pv\_p,$ $\a\_p\ge0,$ $p=1\cdcc m$\\
or total indifference
\ms\ms

}
\\[-.8mm]
\raggedright
Milnor \cite{Miln54}
&&\raggedright
$m$ utility functions $v\_1(x)\cdcc v\_m(x)$, $x\in X$
&&weak order
&&\raggedright
strengthened IIA, N, A, weak P,
independence of the individual zeros
&&{\raggedright
ranking by
$\sum_{p=1}^mv\_p$
\ms\ms\ms

}
\\[-.8mm]
\raggedright
Harsanyi \cite{Hars55}
&&\raggedright
$m$ utility functions $v\_1(x)\cdcc v\_m(x)$, $x\in \Xtil$,\\
$\Xtil$ being the set of lotteries over $X$
&&\raggedright
utility function $v(x)$, $x\in \Xtil$
&&\raggedright
a. Pareto indifference,\\
von~Neumann-Morgenstern utility axioms\\
b. nonstrict Pareto preference + `a.'
&&{\raggedright
a. $v(x)=\sum_{p=1}^m\a\_p\vpx+\beta,$ with some
$\a\_1\cdcc \a\_m,$ $\beta\in\R$\\
b.  `a.' with $\a\_p\ge0,$ $p=1\cdcc m$
\ms\ms\ms

}
\\[-.8mm]
\raggedright
Fishburn \cite{Fish69}
&&$m$ weak orders
&& weak order
&&\raggedright
Pareto principle for equicardinal subsets
&&{\raggedright
ranking by the sum of arbitrary utilities that represent individual
orders

}
\\[-.8mm]
\raggedright
Morkeli\=unas \cite{Mork71}
&&$m$ weak orders
&&weak order
&&\raggedright
N, relaxations of IIA,
weak versions of reinforcement, A, P
&&{\raggedright
a. factored Borda ordering\\
b. extended Borda ordering
\ms

}
\\[-.8mm]
\raggedright G\"ardenfors \cite{Gard73}
&&\raggedright
a. $m$ linear orders\\
b. $m$ weak orders
&&weak order
&&\raggedright
0. N, summability\\
a. weak P, strong positional independence + `0.'\\
b. strong M, stability + `0.'
&&{\raggedright
a. Borda ordering\\
b. extended Borda ordering
\ms\ms\ms\ms

}
\label{TCha}
\end{tabular}\\[-1mm]
\end{Table}
\hyphenpenalty=50
\clearpage

\addtocounter{table}{-1}
\hyphenpenalty=10000
\begin{Table}
{Axiomatic characterizations of scoring methods (Part 2)}
\begin{tabular}[t]{p{1.8cm}|p{.1cm}p{2.19cm}p{.2cm}p{1.64cm}p{.2cm}%
p{3.54cm}p{.2cm}p{3.73cm}}
\begin{tabular}{l}
Paper
\end{tabular}
&&\begin{tabular}{l}
Input\\preferences
\end{tabular}
&&\begin{tabular}{l}
Resulting\\(social)\\structure
\end{tabular}
&&\begin{tabular}{l}
Most important\\axioms
\end{tabular}
&&\begin{tabular}{l}
Result of the method
\end{tabular}
\\[2mm]
\hline
\\[-.8mm]
Smith \cite{Smith73}
&&$m$ linear orders
&&weak order
&&\raggedright
a. reinforcement, N\\
b. overwhelming\\
majority + `a.'
&&{\raggedright
a. composite point ordering\\
b. point ordering

}
\\[-.8mm]
\raggedright
Fishburn \cite{Fish73}
&&\raggedright
$m$ partial orders
&&\raggedright
choice function
&&\raggedright
a. Pareto principle for equicardinal subsets\\
b.,c. versions of A +`a.'\\
d. N + `c.'
&&{\raggedright
a. choice by the sum of arbitrary utilities that weakly represent
individual orders and may depend on the feasible set and the profile\\
b.--d. more symmetric versions of `a.'

}
\\[-.8mm]
\raggedright
Fine, Fine \cite{FinFin74}
&&$m$ weak orders
&&\raggedright
a.,b.,e. weak order\\
c.,d. partial order
&&\raggedright
0. symmetry, elimination, M\\
a. N, contraction + strong `0.'\\
b. overwhelming majority + `a.'\\
c. intersection of all point orderings\\
d. weakest rule to satisfy dissolution + `0.'\\
e. N, inversion, independence of extreme alternatives +
strong `0.'
&&{\raggedright
a. transfinite point ordering\\
b. point ordering\\
c.,d. partial ordering of permuted dominance\\
e. Borda ordering
\ms\ms\ms\ms\ms\ms\ms

}
\\[-.8mm]
Young \cite{Young74}
&&$m$ linear orders
&&choice set
&&\raggedright
reinforcement, N, cancellation, faithfulness
&&{\raggedright
Borda choice
\ms

}
\\[-.8mm]
Young \cite{Young75}
&&$m$ linear orders
&&choice set
&&\raggedright
a. reinforcement, N, A\\
b. overwhelming\\
majority + `a.'
&&{\raggedright
a. composite point choice\\
b. point choice
\ms

}
\\[-.8mm]
\raggedright
Hansson, Sahlquist \cite{HanSah}
&&$m$ linear orders
&&choice set
&&\raggedright
(another proof of Young's \cite{Young74} result)
&&{\raggedright
Borda choice
\ms\ms

}
\\[-.8mm]
Keeney \cite{Keen76}
&&\raggedright
$m$ utility functions $v\_1(x)\cdcc v\_m(x)$, $x\in X$
&&\raggedright
utility function $v(x)$, $x\in X$
&&\raggedright
expected utility assumption,
a version of IIA, M, non-dictatorship
&&{\raggedright
a. $v(x)=\sum_{p=1}^m\a\_p\vpx,$ with $\a\_p\ge0,$ $p=1\cdcc m,$
$\a\_p>0$ for at least two $p$'s
\ms

}
\end{tabular}\\[-1mm]
\end{Table}
\hyphenpenalty=50
\clearpage

\addtocounter{table}{-1}
\hyphenpenalty=10000
\begin{Table}
{Axiomatic characterizations of scoring methods (Part 3)}
\begin{tabular}[t]{p{1.8cm}|p{.1cm}p{2.19cm}p{.2cm}p{1.64cm}p{.2cm}%
p{3.54cm}p{.2cm}p{3.73cm}}
\begin{tabular}{l}
Paper
\end{tabular}
&&\begin{tabular}{l}
Input\\preferences
\end{tabular}
&&\begin{tabular}{l}
Resulting\\(social)\\structure
\end{tabular}
&&\begin{tabular}{l}
Most important\\axioms
\end{tabular}
&&\begin{tabular}{l}
Result of the method
\end{tabular}
\\[2mm]
\hline
\\[-.8mm]
\raggedright
{\dasp} Gevers \cite{dAspGev77}
&&\raggedright
$m$ utility functions $v\_1(x)\cdcc v\_m(x)$, $x\in X$
&&weak order
&&\raggedright
IIA, P, independence of the individual zeros
and common unit, A
&&{\raggedright
ranking by $\sum_{p=1}^m\vpx$ (utilitarianism)
\ms\ms

}
\\[-.8mm]
Maskin \cite{Mask78}
&&\raggedright
$m$ utility functions $v\_1(x)\cdcc v\_m(x)$, $x\in X$
&&weak order
&&\raggedright
a. IIA, P, continuity\\
b. independence of generally unconcerned individuals + `a.'\\
c. A + `b.'\\
d. independence of the common zero and unit of utility + `c.'
&&{\raggedright
ranking by:\\
a. $W(v\_1(x)\cdcc v\_m(x))$ with continuous $W$\\
b. $\sum_{p=1}^mw\_p(\vpx)$ with continuous $w\_1\cdcc w\_m$\\
c. $\sum_{p=1}^mw(\vpx)$ with continuous \& increasing $w$\\
d. $\sum_{p=1}^m\vpx$

}
\\[-.8mm]
\raggedright
Deschamps, Gevers \cite{DesGev78}
&&\raggedright
$m$ utility functions $v\_1(x)\cdcc v\_m(x)$, $x\in X$
&&weak order
&&\raggedright
IIA, P, minimal equity, independence of generally
unconcerned individuals, A,
independence of the common zero and unit
&&{\raggedright
ranking by $\sum_{p=1}^m\vpx$
\ms\ms\ms\ms\ms

}
\\[-.8mm]
\raggedright
Fishburn \cite{Fish78}
&&$m$ choice sets
&&choice set
&&\raggedright
N, reinforcement, disjoint equality
&&{\raggedright
choice set of approval voting (a special case of factored Borda
choice)

}
\\[-.8mm]
\raggedright
Young, Levenglick \cite{YouLev}
&&$m$ linear orders
&&a number of linear orders\raggedright
&&\raggedright
reinforcement, N,
the order of immediately successive alternatives obeys majority vote
&&{\raggedright
Kemeny median
\ms\ms\ms

}
\\[-.8mm]
\raggedright
Richelson \cite{Rich78}
&&$m$ linear orders
&&choice set
&&\raggedright
N, reinforcement, A, independence of Pareto dominated
alternatives
&&{\raggedright
plurality choice
\ms\ms

}
\\[-.8mm]
\raggedright
Farkas, Nitzan \cite{FarNit79}
&&$m$ linear orders
&&choice set
&&\raggedright
closeness to unanimity
&&{\raggedright
Borda choice
\ms

}
\\
\end{tabular}\\[-1mm]
\end{Table}
\hyphenpenalty=50
\clearpage

\addtocounter{table}{-1}
\hyphenpenalty=10000
\begin{Table}
{Axiomatic characterizations of scoring methods (Part 4)}
\begin{tabular}{p{1.8cm}|p{.1cm}p{2.19cm}p{.2cm}p{1.64cm}p{.2cm}%
p{3.19cm}p{.2cm}p{4.08cm}}
\begin{tabular}{l}
Paper
\end{tabular}
&&\begin{tabular}{l}
Input\\preferences
\end{tabular}
&&\begin{tabular}{l}
Resulting\\(social)\\structure
\end{tabular}
&&\begin{tabular}{l}
Most important\\axioms
\end{tabular}
&&\begin{tabular}{l}
Result of the method
\end{tabular}
\\[2mm]
\hline
\\[-.8mm]
Barbera \cite{Barb}
&&$m$ linear orders
&&scores\footnotemark{}
\newcounter{distr}\setcounter{distr}{\value{footnote}}\raggedright
&&\raggedright
a. N, A, strategy- proofness\\
b. separability of individual influences + `a.'\\
c. a relaxed IIA + `a.'\\
d. `b.' + `c.'
&&{\raggedright
a. a convex combination of point scores and lobby size
scores\footnotemark{}\\
\newcounter{comb}\setcounter{comb}{\value{footnote}}
b. point scores$^\arabic{comb}$\\
c. lobby size scores$^\arabic{comb}$\\
d. extended Borda scores$^\arabic{comb}$
\ms

}
\\[-.8mm]
\raggedright
Rubinstein \cite{Rubi}
&&a tournament
&&weak order
&&\raggedright
a relaxed IIA, N, strong M
&&{\raggedright
Copeland ordering
\ms

}
\\[-.8mm]
Roberts \cite{Robe80}
&&\raggedright
$m$ utility functions $v\_1(x)\cdcc v\_m(x)$, $x\in X$
&&weak order
&&\raggedright
0. IIA, weak P\\
a. weak continuity + `0.'\\
b. independence of the common unit of utility + `0.'\\
c. independence of generally unconcerned individuals + `b.'\\
d. independence of the common zero of utility + `b.'\\
e. independence of the individual zeros of utility + `b.'
&&{\raggedright
ranking by:\\
a. $W(v\_1(x)\cdcc v\_m(x))$ with $W$ continuous and increasing\\
b. $W(v\_1(x)\cdcc v\_m(x))$ with $W$ homothetic and increasing\\
c. ${\rm sign}\beta\sum_{p=1}^m\a\_p(\vpx)^\beta$, $\beta\ne0$ or
$\sum_{p=1}^m\a\_p\log(\vpx)$\\
d. $\bar v(x)+W(v\_1(x)-\bar v(x)\cdcc v\_m(x)-\bar v(x)),$ with
$\bar v(x)=1/m\sum_{p=1}^m\vpx$ and $W$ homogeneous of degree 1\\
e. $\sum_{p=1}^m\a\_p\vpx$ with $\a\_p>0,$ $p=1\cdcc m$

}
\\[-.8mm]
\raggedright
Nitzan, Rubinstein \cite{NitRub}
&&$m$ tournaments
&&weak order
&&\raggedright
reinforcement, N, strong M, cancellation
&&{\raggedright
extended Borda ordering
\ms\ms

}
\\[-.8mm]
\raggedright
Myerson \cite{Myer81}
&&\raggedright
a closed, convex and comprehensive\hmii\footnotemark{}
set $\Xtil$ of utility vectors
$\v=(v\_1\cdcc v\_m)$
\newcounter{comp}\setcounter{comp}{\value{footnote}}
&&\raggedright
utility vector $\v^*$ in $\Xtil$
&&\raggedright
weak P, linearity w.r.t.
sets $\Xtil$ (can be replaced by additivity w.r.t.\ sets
\cite{Moul88})
&&{\raggedright
choice of $\v^*=(v^*_1\cdcc v^*_m)\in \Xtil$ to maximize
$\sum_{p=1}^m\a\_pv\_p$ with arbitrary fixed $(\a\_1\cdcc \a\_m)$:
every $\a\_p\ge0$, $\sum_{p=1}^m\a\_p=1$
\ms

}
\\[-1mm]
\hline
\end{tabular}\\
\noalign{\vskip 9pt}
\ \ \ $^\arabic{distr}$Nonnegative scores that sum to 1 (distributions
over $X$)\\
\ \ \ $^\arabic{comb}$The scores are normalized so as to make up
a distribution over $X$\\
\ \ \ $^\arabic{comp}$A set $V\subset\R^m$ is {\em comprehensive\/}
iff $x\in V$ and $y\le x$ (this means $y\_p\le x\_p$ for every\\
$p=\om$) together imply that $y\in V$\\
\noalign{\vskip 6cm}
\end{Table}
\hyphenpenalty=50
\clearpage

\addtocounter{table}{-1}
\hyphenpenalty=10000
\begin{Table}
{Axiomatic characterizations of scoring methods (Part 5)}
\begin{tabular}[t]{p{1.8cm}|p{.1cm}p{2.19cm}p{.2cm}p{1.92cm}p{.2cm}%
p{3.82cm}p{.2cm}p{3.17cm}}
\begin{tabular}{l}
Paper
\end{tabular}
&&\begin{tabular}{l}
Input\\preferences
\end{tabular}
&&\begin{tabular}{l}
Resulting\\(social)\\structure
\end{tabular}
&&\begin{tabular}{l}
Most important\\axioms
\end{tabular}
&&\begin{tabular}{l}
Result of the method
\end{tabular}
\\[2mm]
\hline
\\[-1.0mm]
\raggedright
Morkeli\=unas \cite{Mork82}
&&$m$ weak orders
&&choice set
&&\raggedright
0. reinforcement, N\\
a. independence of\\ Pareto dominated alternatives + `0.'\\
b. independence of the Pareto inferior alternative, duality + `0.'
&&{\raggedright
a. generalized plurality choice\\
b. factored Borda choice
\ms\ms\ms

}
\\[-1.0mm]
Henriet \cite{Henr}
&&\raggedright
a connected relation
&&\raggedright
a. weak order\\
b.,c. choice function
&&\raggedright
0. N, strong M\\
a.,b. independence of cycles + `0.'\\
c. a relaxed IIA + `0.'
&&{\raggedright
a. Copeland ordering\\
b.,c. relative Copeland choice function

}
\\[-1.0mm]
Young \cite{Young86}
&&\raggedright
balanced paired comparisons without ties
&&\raggedright
a. choice set\\
b. a number of linear orders
&&\raggedright
a. N, reinforcement,\\
weak unanimity\\
b. local IIA + `a.'
&&{\raggedright
a. extended Borda choice\\
b. Kemeny median
\ms

}
\\[-1.0mm]
\raggedright
Moulin \cite{Moul88}
&&\raggedright
a closed, convex and comprehensive\hmii\footnotemark{}
set $\Xtil\subset \R^m$ of utility
vectors $\v=(v\_1\cdcc v\_m)$
\newcounter{compm}\setcounter{compm}{\value{footnote}}
&&\raggedright
utility vector $\v^*$ in $\Xtil$
&&\raggedright
NIIA, P, A, commutativity with translations
&&{\raggedright
choice of $\v^*=(v^*_1\cdcc v^*_m)\in \Xtil$ to maximize
$\sum_{p=1}^mv\_p$
\ms\ms\ms

}
\\[-1.0mm]
\raggedright
Chebotarev \cite{Cheb89a,Cheb94}
&&\raggedright$m$ skew-symmetric incomplete matrices of paired comparisons
&&scores
&&\raggedright
coincidence with extended Borda scores for complete paired
comparisons,\\
implicit form of scores
&&{\raggedright
generalized row sums
\ms\ms\ms\ms\ms

}
\\[-1.0mm]
Saari \cite{SaariSCW}
&&$m$ linear orders
&&choice set
&&\raggedright
relaxed versions of Young's \cite{Young74} axioms
&&{\raggedright
Borda choice
\ms

}
\\[-1.0mm]
\raggedright
Barzilai, Golany \cite{BarGol}
&&\raggedright
a skew- symmetric matrix $(a\_{ij})$ of paired comparisons
&&scores
&&\raggedright
N, additivity, recovery of ${\bf s}\in\R^n$ such that
$a\_{ij}=s\_i-s\_j,$ $i,j=1\cdcc n$ (whenever ${\bf s}$ exists)

&&{\raggedright
normalized extended Borda scores
\ms\ms\ms

}
\\[-1.0mm]
\raggedright
Baigent, Xu \cite{BaiXu91}
&&\raggedright
$m$ choice functions
&&\raggedright
choice function
&&\raggedright
N, M, independence of symmetric substitutions
&&{\raggedright
choice function of approval voting

}
\\[-1mm]
\hline
\end{tabular}\\
\noalign{\vskip 7pt}
\ \ \ $^\arabic{compm}$See footnote \arabic{comp}\\
\noalign{\vskip 8cm}
\end{Table}
\hyphenpenalty=50
\clearpage

\addtocounter{table}{-1}
\hyphenpenalty=10000
\begin{Table}
{Axiomatic characterizations of scoring methods (Part 6)}
\begin{tabular}[t]{p{1.8cm}|p{.1cm}p{2.19cm}p{.2cm}p{1.64cm}p{.2cm}%
p{3.54cm}p{.2cm}p{3.73cm}}
\begin{tabular}{l}
Paper
\end{tabular}
&&\begin{tabular}{l}
Input\\preferences
\end{tabular}
&&\begin{tabular}{l}
Resulting\\(social)\\structure
\end{tabular}
&&\begin{tabular}{l}
Most important\\axioms
\end{tabular}
&&\begin{tabular}{l}
Result of the method
\end{tabular}
\\[2mm]
\hline
\\[-.2mm]
Debord \cite{Debo}
&&$m$ binary relations\footnotemark{}
\newcounter{Deb}\setcounter{Deb}{\value{footnote}}\raggedright
&&\raggedright
a number of $k$-choice sets
&&\raggedright
reinforcement, N, cancellation, faithfulness
&&{\raggedright
$k$-sets with maximal total extended Borda score
\ms

}
\\[-.2mm]
\raggedright
Bouyssou\\
\cite{Bouy92a}
&&\raggedright
a valued relation
&&weak order
&&\raggedright
N, strong M, independence of circuits
&&{\raggedright
extended Borda ordering
\ms

}
\\[-.2mm]
\raggedright
Bouyssou\\
\cite{Bouy92b}
&&\raggedright
a valued relation
&&\raggedright
choice function
&&\raggedright
N, strong M, independence of circuits
&&{\raggedright
relative extended Borda choice function

}
\\[-.2mm]
\raggedright
Bouyssou, Perny \cite{BouPer}
&&\raggedright
a valued relation
&&\raggedright
partial order
&&\raggedright
relaxed N, strong M, independence of alternated cycles
&&{\raggedright
meet of down-sided and up-sided Borda orderings
\ms

}
\\[-.2mm]
Duggan \cite{Dugg}
&&\raggedright
$m$ weak orders
&&\raggedright
choice function
&&\raggedright
N, weak P, A,
a rank version of independence of the individual zeros
&&{\raggedright
absolute down-sided Borda choice function
\ms

}
\\[-.2mm]
\raggedright
Myerson\\
\cite{Myer}
&&\raggedright
$m$ arbitrary ballots
&&choice set
&&\raggedright
a relaxed N, A, reinforcement,
overwhelming majority
&&{\raggedright
choice of alternatives maximizing the sum of arbitrary
scores determined by individuals' ballots

}
\\[-.2mm]
\raggedright
Marchant \cite{Marc96a}
&&\raggedright
$m$ valued relations\footnotemark{}
\newcounter{Mar}\setcounter{Mar}{\value{footnote}}
&&weak order
&&\raggedright
N, reinforcement, cancellation, faithfulness
&&{\raggedright
extended Borda ordering
\ms

}
\\[-.2mm]
\raggedright
Marchant \cite{Marc96b}
&&\raggedright
a. $m$ valued relations$^\arabic{Mar}$\\
b. $m$ rational valued relations$^\arabic{Mar}$\\
&&scores
&&\raggedright
0. N, cancellation, faithfulness\\
a. reinforcement + `0.'\\
b. a relaxed reinforcement + `0.'
&&{\raggedright
a.,b. extended Borda scores up to a positive affine transformation

}
\\[-1mm]
\hline
\end{tabular}\\
\noalign{\vskip 9pt}
\ \ \ $^\arabic{Deb}$The set of admissible relations
contains all linear orders and is stable to transpositions\\
\ \ \ $^\arabic{Mar}$The values are from $[0,1]$ (fuzzy relations).
The set of admissible relations contains all\\
crisp weak orders and is stable to transpositions and permutations\\
\noalign{\vskip 8cm}
\end{Table}
\hyphenpenalty=50
\clearpage

\subsection{Input preferences}

The form of input information also varies:

\begin{itemize}
\item Typically, it is a profile of classical individual preferences:
     \begin{itemize}
     \item Linear orders
          \cite{Gard73,Smith73,Young74,Young75,HanSah,YouLev,Rich78,%
FarNit79,SaariSCW}.
     \item Weak orders
        \cite{Flem52,Fish69,Mork71,Gard73,FinFin74,Mork82,Dugg}.
     \item Individuals' scores, points, marks, cardinal utilities
        \cite{GooMar,BlaGir54,Miln54,Hars55,Keen76,dAspGev77,Mask78,%
DesGev78,Robe80,Myer81,Moul88}.
     \item Choice sets \cite{Fish78}.
     \end{itemize}
\item The preferences in the profile may have more general forms:
     \begin{itemize}
      \item Partial orders \cite{Fish73}.
      \item Tournaments \cite{NitRub}.
     \item Arbitrary binary relations \cite{Debo} (see also
                                footnote~\arabic{Deb}). 
     \item Valued (fuzzy) relations \cite{Marc96a,Marc96b} (see
                                footnote~\arabic{Mar}). 
     \item Incomplete skew-symmetric matrices of paired comparisons (in
     other words, skew-symmetric valued relations with incomparability
     distinguished from zero values) \cite{Cheb89a,Cheb94}.
     \item Choice functions \cite{BaiXu91}.
     \item Arbitrary individual ballots \cite{Myer}.
     \end{itemize}
\item In some papers, the initial preferences are represented by a
      single relation (which may be thought of as a majority relation
      or other function of the individual profile):
     \begin{itemize}
     \item Tournament \cite{Rubi}.
     \item Connected relation \cite{Henr}.
     \item Valued relation \cite{Bouy92a,Bouy92b,BouPer}.
     \item Skew-symmetric matrix of paired comparisons \cite{BarGol}.
     \end{itemize}
\item Young \cite{Young86} aggregates a series of paired comparisons
     without ties, where every alternative is involved in the same
     number of comparisons.
\end{itemize}

\subsection{The most important axioms}

Many of the following axioms are applied in slightly varying versions
in different papers. Therefore we prefer to give their main ideas
rather than exact formulations. Some axioms are renamed for the sake of
unification.

\bigskip{\bf Variable electorate axioms}.
{\em Reinforcement\/} for choice procedures states that if $\A$ and
$\B$ are disjoint preference profiles, $\A+\B$ is the combined profile
and $C(\A)$ is the choice set for $\A$, then $C(\A+\B)=C(\A)\bigcap
C(\B)$ whenever $C(\A)\bigcap C(\B)\ne\emptyset$. The same condition
expresses reinforcement for procedures that generate a number of
orders; here $C(\A)$ denotes the set of these orders. For ranking
procedures, reinforcement states that if $i$ is socially no worse than
$j$ for both $\A$ and $\B$, then this is the case for $\A+\B$ (with
strict social preference for $\A$ or $\B$ implying strictness for
$\A+\B$).  Reinforcement for scoring procedures means that every
alternative's score for $\A+\B$ is the sum of its scores for $\A$ and
$\B$, i.e., it reduces to a kind of additivity. Obviously, this axiom
coincides with the following one if they both are applied to matrices
of summarized paired comparisons.

{\em Additivity\/} for scoring procedures that operate on paired
comparison matrices means that the sum of two matrices is mapped to the
sum of score vectors derived from these matrices.

The following two axioms are `dual' to reinforcement.

{\em Elimination\/} says that if $i$ is socially no worse than $j$ for
$\A+\B$ and $i$ is socially equivalent to $j$ for $\A$, then $i$ is
socially no worse than $j$ for $\B$.

{\em Strong elimination\/} says that if $i$ is socially no worse than
$j$ for $\A+\B$ and $j$ is socially no worse than $i$ for $\A$, then
$i$ is socially no worse than $j$ for $\B$.

{\em Overwhelming majority\/} for choice procedures states that whenever
$C(\A)=\{i\}$, then for any $\B$ there is an integer $k^*$ such that
$C(k\A+\B)=\{i\}$ for all $k\ge k^*$. Here $k\A$ is the profile
consisting of $k$ copies of $\A$. Overwhelming
majority for ranking procedures results if we replace choice by social
binary preferences.

{\em Contraction\/}, in the case of procedures that produce weak
orders, says that the social ordering for $\A$ coincides with that for
$k\A$.

{\it Disjoint equality\/} says that if there are only two individuals,
their individual choice sets $X^1$ and $X^2$ are nonempty (which is not
generally assumed), and $X^1\cap X^2=\emptyset$, then $X^1\cup X^2$ is
exactly the social choice. This axiom is important for characterizing
approval voting \cite{Fish78}. Another related axiomatization is given
in \cite{Fish79}.

{\em Faithfulness\/} states that in the case of only one individual
having a linear order, her top ranked alternative/alternatives
constitute the social choice/$k$-choice (respectively, individual's
weak order is taken as the social one or is concordant with the social
scores).

This natural condition is a variable electorate axiom of positive
relation between individual and social preferences. Other axioms of
this kind follow.

\bigskip{\bf Axioms of positive relation between individual and
social preferences}.
Monotonicity (`M' in Table~\ref{TCha}) says that if some alternative
becomes more favorite in one individual opinion, whereas all other
alternatives get no rise, then it does not become worse in the social
preference.  Slightly varying formulations of {\em strong
monotonicity\/} (`strong~M' in Table~\ref{TCha}) additionally require
that this alternative leave behind the alternatives that were socially
indifferent to it.

The following axiom restricts the positive response of the social
preference.

Suppose that some individual $p$ only changes $x\prec\_p y$ to $x\sim\_p
y$ or $x\sim\_p y$ to $x\succ\_p y$ and others change nothing.
{\em Stability\/} applies to the procedures that transform profiles of
weak orders to social weak orders and says that in no above situation
the social weak order can change $x\prec z$ to $x\succ z$ with any
$z\ne y$.

{\em Pareto indifference\/} says that everyone being indifferent
between two alternatives implies social indifference between them.

{\em Pareto preference\/} says that if no individual strictly prefers
$j$ to $i$ and at least one strictly prefers $i$ to $j$, then $i$ is
strictly socially preferred to $j$.

{\em Nonstrict Pareto preference\/} says that if no individual prefers
$j$ to $i$, then the society cannot strictly prefer $j$ to $i$.

{\em Strong Pareto principle\/} (`P' in Table~\ref{TCha}) is the
conjunction of Pareto indifference and Pareto preference.

{\em Weak Pareto principle\/} says that if all individuals strictly
prefer $i$ to $j$ then so does the society.

To introduce the following axiom, suppose $\succeq$ is a weak order on
the set of alternatives $X$ and $\sim$ is its indifference part. We say
that a subset of alternatives $X^1\subseteq X$ is superior (equivalent)
to $X^2\subseteq X$ w.r.t.\ $\succeq$ iff ${\rm card}X^1={\rm card}X^2$
and there exists a one-to-one correspondence $\omega$ from $X^1$ onto
$X^2$ such that $\omega(i\_1)=i\_2$ implies $i\_1\succeq i\_2$ (resp.,
$i\_1\sim i\_2$). Strict superiority means superiority and not
equivalence. {\em Pareto principle for equicardinal subsets\/} says
that whenever $X^1$ is superior to $X^2$ w.r.t.\ every individual's
weak order, then $X^2$ cannot be strictly superior to $X^1$ w.r.t.\ the
social weak order. If, in addition, $X^1$ is not equivalent to $X^2$
w.r.t.\ at least one individual weak order, then $X^2$ cannot be
equivalent to $X^1$ w.r.t.\ the social weak order.

{\em Weak unanimity\/} for choice procedures says that whenever all
individual paired comparisons favor alternative $i$, then $\{i\}$ is
the social choice set. Weak unanimity for ranking procedures says that
whenever all individual paired comparisons agree with a fixed linear
order, then it is taken as the social one.

{\em Closeness to unanimity\/} essentially specifies a concrete voting
procedure which follows. Let $\U(i)$ be the set of all profiles where
$i$ is the top alternative in all constituent linear orders. For a
given profile $\A$, those alternatives $i$ are chosen which minimize
the {\em inversion distance} between $\A$ and $\U(i)$.  As usual, this
distance is defined as that between $\A$ and the nearest member of
$\U(i)$ and equals the number of differently ordered pairs of
alternatives.

Table~\ref{TCha} presents only one characterization \cite{FarNit79}
based on resolving optimization problems. In fact, these are numerous.
One of the most interesting results of this nature is as follows.

The Borda method can be characterized \cite{Newe92} as the point method
which maximizes (among all point methods) the proportion of profiles
for which the social preferences agree with those of the majority vote
(viz., with the majority winner on a pair of alternatives, or with the
Condorcet winner, or with the Condorcet ordering, provided that the
latter two exist). A related characterization of the Borda choice
function was obtained by Fishburn and Gehrlein \cite{FisGer76}. Some
other characterizations involving optimization are dealt with in
\cite{CooSei82,Cheb89b,Cheb90a,CheSha97}; \cite{Cheb89b} also describes
some statistical procedures resulting in Borda-like scores.

\bigskip{\bf Permutation-independence and substitution-independence
axioms}.
{\em Ano\-nymity\/} and {\em neutrality\/} essentially say that no
information about individuals and alternatives (respectively) except
for the preference profile is used to derive the social preferences.
More formally, they require the social preference operator to be
stable to any permutation of individuals and commutative with any
permutation of alternatives.

Suppose that there exist a permutation of individuals and a permutation
of alternatives such that their simultaneous application leaves the
preference profile invariant. {\em Symmetry\/} states that the
corresponding alternatives must be socially equivalent in this case.

Symmetry is very effective as applied to the voting situations like the
Condorcet paradox. On the relationship between anonymity, neutrality
and symmetry see \cite[pp.~473--474]{FinFin74}.

{\it Independence of symmetric substitutions says\/}: If profile $\B$
is obtained from $\A$ by a one-element swap between $C^p(X')$ and
$C^q(X')$ (individual choices of $p$ and $q$ from feasible set
$X'\subseteq X$), then the social choice from $X'$ must be the same
for $\A$ and $\B$.

\bigskip{\bf Individual-independence axioms}.
{\em Independence of locally unaffected individuals\/} \cite{Flem52}
states that if some individual is indifferent between two alternatives,
then his/her other preferences do not influence the social preference
between these two alternatives.

{\em Independence of generally unconcerned individuals\/} says that if
some individual assigns the same score (utility) to all alternatives,
then the social ordering is not affected by this particular score. This
corresponds to Debreu's \cite{Debr60} strong separability condition
which is crucial for his derivation of the additive utility
representation.

Independence of locally unaffected individuals plays a similar role in
Fleming's derivation of monotone summability.

Suppose that only one individual changes his/her preference relation.
{\em Separability of individual influences\/} applies to scoring
procedures and says that the resulting differences of the scores of
alternatives solely depend on that individual's previous and new
preferences and are indifferent to the other individuals' preferences.

\bigskip{\bf Alternative-independence axioms}.
{\em Arrow's independence of irrelevant alternatives\/} (`IIA' in
Table~\ref{TCha}) says that the application of the social preference to
any subset of the set of alternatives solely depends on the restriction
of the individual profile to this subset.

{\em Strengthened IIA\/} says that after the removal of any
alternative, the social preference coincides with the initial one
applied to the subset of remaining alternatives, provided that the
individual preferences do not change.

As distinct from IIA, strengthened IIA implies that the set (and the
number) of alternatives may vary.

{\em Nash's independence of irrelevant alternatives\/} (NIIA in
Table~\ref{TCha}) says that whenever $X'$ and $X''\subset X'$ are two
sets of alternatives such that $X''$ contains some alternatives chosen
from $X'$, then exactly these alternatives constitute the choice from
$X''$.

Nash applied this condition to one-element choice, and the more general
above formulation is due to Arrow.

Note that IIA is met by (and enters characterizations of) most scoring
methods that operate on individual utilities but is not satisfied by
many scoring methods dealing with binary relations (including valued
relations). This is due to the very {\em relative nature\/} of binary
relations which causes some loss of information involving a subset
of alternatives as a profile is restricted to the subset. The following
observation clarifies the point. Consider the restrictions of an
individual profile to two complementary subsets of alternatives. Then
the initial profile can be retrieved from these two restrictions
provided that it consists of utilities and cannot be if it comprises
binary relations. From this point, the variable electorate axioms are
more appropriate in the case of binary relations, since no information
is lost when the profile is divided into parts corresponding to
disjoint sets of individuals. The same can be said of the profiles that
consist of choice functions.

The record `a relaxed IIA' stands for various relaxations of IIA.
The following two axioms combine relaxed IIA with a positive relation
to social preferences.

{\it Independence of Pareto dominated alternatives\/} says that
removing Pareto dominated alternatives does not alter the social
choice.

{\em Independence of the Pareto inferior alternative\/} states the same
concerning the alternative that is Pareto dominated by every other one.

{\em Independence of extreme alternatives\/} is a close condition
which, unlike the previous two axioms, does not imply any `positive
relation to social preferences'. It says that unanimously superior as
well as unanimously inferior subsets of alternatives can be discarded
without changing the order on the remaining alternatives.

The following axiom embodies a similar idea applied to one individual's
preferences.

Suppose that some individual ranks alternative $z$ ahead both $x$ and
$y$ or behind both $x$ and $y$, and then moves $z$ to another position
which is also not between $x$ and $y$ and not the same with $x$ or $y$.
{\em Strong positional independence\/} for the profiles of weak orders
\cite{Gard73} says that such a shift does not change the social
relation between $x$ and $y$.

{\em Dissolution\/} \cite{FinFin74} introduces acceptable dissolution
procedures (we do not specify them here) for breaking ties in the
individual profiles of weak orders and says that the social binary
relation for a given profile must contain the common part of social
binary relations derived for all acceptable dissolutions of the
profile.

A similar idea has been exploited in \cite{Mork71}.

Young and Levenglick \cite{YouLev} use a peculiar axiom which
combines independence and positive reaction. In Table~\ref{TCha} it
reads as {\em `the order of immediately successive alternatives obeys
majority vote'.}  This regards every two alternatives immediately
successive in one of the derived social linear orders.

\bigskip{\bf Independence axioms involving preference relativity}.
{\em Independence of the individual zeros (of utility)\/} requires that
the social ordering over $X$ remain unchanged when each individual's
utility function $v\_p(x),$ $x\in X,$ $p=\om$ is replaced by
$v'_p(x)=v\_p(x)+c\_p$ with any $c\_1\cdc c\_m$.

{\em Independence of the common zero\/} and {\em independence of the
common unit (of utility)\/} respectively mean the invariance of the
social ordering with respect to arbitrary transformations
$v'_p(x)=v\_p(x)+c,$ $\;c\in\R$ and $v'_p(x)=bp(x),$ $\;b>0$,
$\;p=\om$, $\;x\in X$.

{\em Independence of cycles\/} states that the choice set
(respectively, the resulting ordering) does not alter whenever any cycle
in any input preference relation is reversed. {\em Independence of
circuits\/} and {\em independence of alternated cycles\/} are
variations of this axiom.

The following condition is closely related to independence axioms. In
the case of ordinary binary relations it can be termed `majority
equivalence', whereas in the case of valued relations it bears the
spirit of additivity.

{\em Cancellation\/} for procedures that operate on profiles of binary
preferences states that if for every $i$ and $j$, sum of the entries of
the relations on $(i,j)$ is equal to that on $(j,i)$, then all
alternatives are socially equivalent. For ordinary binary relations,
`sum of the relations on $(i,j)$' reduces to the number of individuals
that prefer $i$ to $j$.

\bigskip{\bf Other profile transformation axioms}.
Note that monotonicity and independence conditions are profile
transformation axioms. Here we present such transformation conditions
that bear neither positive response nor independence.

{\em Duality\/} states that the reversal of all individual preference
relations cannot preserve any formerly chosen alternative in the social
choice set, unless the initial social choice set coincides with $X$.

{\em Inversion\/} says that the social weak order derived after the
above transformation is the reversed initial one.

Suppose that the alternatives are identified with the vectors of
individual utilities attached to them. Then every set in $\R^m$ can be
considered as a set of alternatives. Let $\vp$ be a social choice
procedure which indicates one `best' alternative (point in $\R^m$) in
each set that belongs to its domain.

{\em Linearity with respect to sets\/} means that for every
$V',V''\in\R^m$, if $V'$ and $V''$ belong to the domain of a social
choice procedure $\vp$, then for any $0\le\lambda\le1,$ $\lambda V'+
(1-\lambda)V''$ also belongs to the domain of $\vp$ and $\vp(\lambda
V'+ (1-\lambda)V'')=\lambda\vp(V')+ (1-\lambda)\vp(V'')$. Here $\lambda
V'+(1-\lambda)V''$ designates $\{\lambda \v'+(1-\lambda)\v''\mid \v'\in
V',\v''\in V''\}$.

This axiom can be given a probabilistic interpretation. If tomorrow one
will have to choose from $V'$ or $V''$ with probabilities $\lambda$ and
$1-\lambda$, respectively, then within the von Neumann-Morgenstern
framework, $\lambda V'+(1-\lambda)V''$ is the set of all expected
utility allocations which are now feasible. Thus, linearity w.r.t.\
sets says that the today's chosen utility allocation should be equal to
the expectation of the tomorrow's allocation.

{\em Additivity with respect to sets\/} is a close property saying that
for any admissible $V'$ and $V''$, $V'+V''$ (i.e., $1\cdot V'+1\cdot
V''$) is also admissible and $\vp(V'+V'')=\vp(V')+\vp(V'')$.

{\em Commutativity with translations\/} reads as the application of
additivity w.r.t.\ sets when one of the sets is a singleton:
$\vp(V+\{v\})=\vp(V)+v$ for all admissible $V$ and all $v\in\R^m$.

\bigskip{\bf Expected utility axioms}.
As Myerson \cite{Myer81} demonstrated and as we saw above, linearity
with respect to sets has a natural interpretation in terms of expected
utility.

In general, if an aggregation procedure takes each profile of expected
utility functions to an expected utility function, the linearity of
this procedure follows under rather weak conditions (namely,
von~Neumann-Morgenstern utility axioms and Pareto indifference are
sufficient).  Harsanyi \cite{Hars55} was the first to show this (see
also \cite{Patt68,Keen76,Fish84}). An excellent philosophical
justification of the underlying model is given in \cite{Hars75,Hars82}.
The entire book \cite{SenWil82} is very valuable for the comprehension
of utilitarianism.

The von~Neumann-Morgenstern utility axioms \cite{vonNeuMor47} can
be summarized (in the form of Marschak postulates \cite{Marsc50}) as
follows. The relation of preference is a weak order ($\succeq$) defined
on a set of alternatives (prospects); this set is closed w.r.t. taking
lotteries (probabilistic mixtures) of its members; there are at least
four mutually nonindifferent prospects; if $x\succeq y\succeq z$ then
there exists a mixture of $x$ and $z$ such as to be exactly indifferent
to $y$; if $x$ and $x'$ are indifferent, then, for any prospect $y$,
any mixture of $x$ and $y$ is indifferent to the mixture of $x'$ and
$y$ with the same probabilities.

The {\em Expected utility assumption\/} of \cite{Keen76} states that
both the individual scores of alternatives and the social scores are
expected utilities subject to the same underlying probabilistic model.
Essentially, this means that the feasible alternatives are lotteries
over certain outcomes.

In fact, Harsanyi's theorem (as well as more recent Keeney's theorem
\cite{Keen76} and the result in \cite{Fish84}) can be considered as a
realization of the program outlined by Fleming \cite[p.~380]{Flem52} on
translating his own result into the expected utility framework. On the
connections between these two results see \cite{Flem57}.

{\em Strategy-proofness\/} means in \cite{Barb} that if the social
scores (properly normalized) are considered as the probabilities of the
alternatives in a lottery, and each individual has a cardinal utility
function that induces her personal ordering, then no individual can
increase her expected utility in the resulting lottery by
misrepresenting her true ordering.

\bigskip{\bf Other axioms}.
Suppose that the individual profile consists of utility functions
$v\_1(i)\cdc v\_m(i)$ on $X$, $\v=(v\_1(i)\cdc v\_m(i))$ is a {\em
utility mapping}, $V$ is the set of admissible utility mappings, and
$\vp$ is a procedure that takes each $\v\in V$ to an ordering $R$ over
$X$.

Let $\x=(x\_1\cdc x\_m)\in\R^m$. Suppose $R^+(\x)$ is the set of
$\y=(y\_1\cdc y\_m)\in\R^m$ such that for some admissible utility
mapping $\v\in V$, $(x\_1\cdc x\_m)$ and $(y\_1\cdc y\_m)$ are the
vectors of individual utilities of some alternatives $i\in X$ and
$j\in X$, respectively, and $iRj$, where $R=\vp(\v)$. The set
$R^-(\x)$ is defined by replacing `$iRj$' with `$jRi$' in the
definition of $R^+(\x)$. {\em Continuity\/} states that the sets
$R^+(\x)$ and $R^-(\x)$ are closed for all $\x\in \R^m$.

It follows from Debreu's teorem \cite{Debr54} that this axiom implies
the representability of $\vp$ by a continuous function of $v\_1(i)\cdc
v\_m(i)$.

Roberts \cite{Robe80} exploits another condition,
{\em weak continuity\/}, which says that for every admissible utility
mapping $\v=(v\_1(i)\cdc$ $v\_m(i))$, there exists an admissible utility
mapping $\v'=(v'_1(i)\cdc$ $v'_m(i))$ that is close to $v$ as much as
desired and such that $v'_p(i)<v\_p(i)$ for all $i$ and $p$ and
$\vp(\v')=\vp(\v)$.

A scoring procedure satisfies {\em summability\/} if there exists
a function $w$ such that the scores of alternatives can be represented
as
$$
s\_i=\sum_{p=1}^m w(A^{(p)},i), \quad i=\on.
$$
Summability of a ranking (choice) procedure means that the social
order (choice) is determined by a summable scoring procedure.

This condition is the crucial axiom in \cite{Gard73} and simultaneously
is closely related to the procedures characterized in
\cite{Flem52,Fish69,Mask78} and point ranking procedures characterized
in \cite{Smith73,FinFin74,Young75,Barb} and generalized in \cite{Myer}.

The idea of `minimal equity' is that the social ranking procedure is not
always based on the most optimistic estimates of alternatives (in other
terms, on the preferences of the best off individual). Specifically,
suppose $i,j\in X$ and there exists $p\in M$ such that $\vpj>\vpi$
whereas for all other $q\in M$, $v\_q(j)<v\_q(i)<v\_p(i)$.  {\em Minimal
equity\/} says that the social ranking procedure is not one always
preferring $j$ in such situations.
\bigskip

Thus, the characterization results rest upon:
\begin{itemize}
\item Reinforcement or elimination (sometimes supplemented by
      overwhelming majority)
      \cite{Smith73,FinFin74,Young74,Young75,HanSah,Fish78,Rich78,%
NitRub,YouLev,Mork82,Young86,SaariSCW,Debo,Myer,Marc96a,Marc96b}.
\item Pareto principle for equicardinal subsets \cite{Fish69,Fish73}.
\item Closeness to unanimity \cite{FarNit79}.
\item Independence of cycles and its variations
      \cite{Henr,Bouy92a,Bouy92b,BouPer}.
\item Independence of locally unaffected individuals \cite{Flem52}.
\item Relaxations of IIA \cite{Rubi,Henr}. Sometimes these conditions
      combine with sum\-mability \cite{Gard73}, elimination
      \cite{FinFin74}, reinforcement \cite{YouLev,Rich78,Young86} and
      its weak version \cite{Mork71}.
\item Additivity \cite{BarGol}, additivity w.r.t.\ sets \cite{Moul88},
      linearity w.r.t.\ sets \cite{Myer81}, commutativity with translations
      \cite{Moul88}.
\item IIA applied to utility functions along with scale invariance
      conditions \cite{GooMar,BlaGir54,Miln54,dAspGev77,Mask78,%
DesGev78,Robe80}.
      In a discrete context, this approach is applied in
      \cite{GooMar,Dugg}.
\item Expected utility axioms \cite{Hars55,Keen76,Fish84}.
\item An expected utility version of strategy-proofness \cite{Barb}.
\item Among other axioms we mention
     \begin{itemize}
     \item Axioms of positive relation between individual and social
           preferences (Pareto principle, monotonicity, independence of
           Pareto dominated/ Pareto inferior alternatives, weak
           unanimity, faithfulness) which occur in almost all
           characterizations.
     \item independence of generally unconcerned individuals
           \cite{Mask78,DesGev78,Robe80} which, as well as independence
           of locally unaffected individuals and Pareto indifference
           (to some extent), can help separate individual variables.
     \item Cancellation \cite{Young74,NitRub,Debo,Marc96a,Marc96b}.
     \item Disjoint equality \cite{Fish78}.
     \item Independence of symmetric substitutions \cite{BaiXu91}.
     \item Symmetry, contraction, dissolution \cite{FinFin74}.
     \item Inversion \cite{FinFin74} and duality \cite{Mork82}.
     \item Summability and stability \cite{Gard73}.
     \item Continuity \cite{Mask78} and weak continuity \cite{Robe80}.
     \item Minimal equity \cite{DesGev78}.
     \item Implicit form of scores \cite{Cheb89a,Cheb94} whose
     generalization will be discussed in Section~\ref{MIF}.
     \end{itemize}
\end{itemize}
\medskip

\subsection{Resulting methods}

The last column of Table~\ref{TCha} requires some definitions. We
represent individual binary relations (both ordinary and valued ones)
by matrices of paired comparisons $A^{(p)}=(a_{ij}^p)$, $p=\on$. When
the input preference relation is single, one matrix $A=(a\_{ij})$ can
be substituted for $A^{(p)}$. As distinct from all other sections of
this paper, now we do not assume
$a_{ij}^p+a_{ji}^p=1$ for $j\ne i$. For ordinary relations $R^p$,
$$
a_{ij}^p=\cases{1, &$(i,j)\in R^p$ and $(j,i)\notin R^p$\cr
                0, &otherwise.\cr}
$$

First, we introduce a rather general form of scores used in the methods
of Table~\ref{TCha}. {\em Extended Borda scores\/} are defined as
follows:
\begin{equation}
\si=\suml_{j\ne i}\suml_{p=1}^m(a_{ij}^p-a_{ji}^p), \qquad i=\on.
\label{EBS}
\end{equation}

The idea of summing differences $(a_{ij}^p-a_{ji}^p)$ is due to
Copeland \cite{Cope} (see also \cite{LucRai}), but conventionally his
name is only attributed to such scores derived from a single relation,
particularly from that of simple majority. We follow this tradition.

{\em Extended Borda ordering\/} and {\em extended Borda choice\/}
designate in Table~\ref{TCha} various instances of ranking and choice
induced by the extended Borda scores, provided that they have no proper
name. Specifically, {\em Borda choice\/} and {\em Borda ordering\/} are
the instances where the individual profile consists of linear orders.
When a choice function is concerned, the Borda (or Copeland) choice may
have {\em absolute\/} (`broad') or {\em relative\/} (`narrow') form
\cite{Sen77a}, depending on whether the score is counted over the whole
set of alternatives or over its presented subset.

The scores
$$
\si=\suml_{j\ne i}\suml_{p=1}^m a_{ij}^p, \qquad i=\on
$$
and
$$
\si=\suml_{j\ne i}\suml_{p=1}^m(-a_{ji}^p), \qquad i=\on
$$
are referred to as {\em down-sided\/} and {\em up-sided Borda scores},
respectively. For profiles of linear orders, they are equivalent to the
extended Borda scores up to a positive affine transformation. In this
case, the term {\em `Borda score'\/} is used. {\em Factored Borda
scores} are also equivalent to the extended Borda scores on all
profiles of linear orders. In the case of weak orders, they are defined
as follows:
$$
s\_i=\sum_{p=1}^ms_i^p, \qquad i=\on,
$$
where
$$
s\_i^p={\rm card}\{k\mid\exists i\_1\cdc i\_k\in X:
i\succ\_pi\_1\succ\_p\cdots\succ\_pi\_k\}
$$
and $\succ\_p$ is the strict preference of the $p$th individual. The
properties of these and some other extensions of the Borda scores to
the profiles of weak orders are studied in \cite{Gard73}.

Suppose that the preference profile consists of individuals' choice
sets.  Note that they can be represented by weak orders with exactly
two `strata' (classes of equivalent alternatives). In this case, the
factored Borda scores are {\em approval voting scores}. They equal the
numbers of supporting individuals.

For profiles consisting of linear orders, down-sided Borda scores are
generalized in two ways. The scores
$$
\si=\suml_{p=1}^m\a(\suml_{j\ne i} a_{ij}^p), \qquad i=\on
$$
with $\a(\cdot)$ nonnegative and nondecreasing\footnote{%
Note that in some characterizations these properties are not derived}
real-valued function, are {\em point scores}.  Function $\a(\cdot)$
defines the `points' assigned to an alternative $i$ for each position
in individual orders. If each position except for the highest one is
assigned a zero point, then {\em plurality scores\/} result. {\em
Transfinite point ordering\/} allows transfinite `points' which
need not satisfy the Archimedian property.  {\em Composite point
ordering\/} results when one has several different vectors of point
scores (determined by different $\a$'s) and applies them successively,
for breaking ties that survived. In fact, transfinite point orderings
and composite point orderings give rise to the same family of
aggregation operators (see \cite{Smith73}). At {\em composite point
choice}, the point scores are used to refine the choice. As distinct
from the point ordering, {\em ranking by the sum of arbitrary scores
(utilities) that represent individual orders\/} is generally neither
anonymous nor neutral. A utility function $v: X\to \R$ {\em represents
a weak order\/} $\succeq$ on $X$ iff $\forall i,j\in X$, $i\succeq j
\equi v(i)\ge v(j)$.

The scores
$$
\si=\suml_{j\ne i}\beta(\suml_{p=1}^m  a_{ij}^p), \qquad i=\on
$$
with $\beta(\cdot)$ nonnegative and nondecreasing real-valued function,
are {\em lobby size scores} ({\em supporting size scores\/} in
\cite{Barb}).  Function $\beta(\cdot)$ defines the partial scores
assigned to $i$ for each size of a `lobby' (coalition) supporting $i$
against $j$. This is another generalization of down-sided Borda scores.

The most general Borda-like scores for profiles of linear orders are
provided by {\em convex combinations\/} of point scores and lobby size
scores:
$$
\si=   \nu \suml_{p=1}^m    \a(\suml_{j\ne i} a_{ij}^p)+
     (1-\nu)\suml_{j\ne i}\beta(\suml_{p=1}^m  a_{ij}^p), \qquad i=\on,
$$
where $0\le\nu\le1$. The same approach can be applied to the extended
Borda scores as well.

{\em Partial ordering of permuted dominance\/} is an interesting
procedure introduced in \cite{FinFin74} as follows. Alternative $i$
nonstrictly dominates $j$ iff there exists a permutation $\sigma$ of
individuals and, for each individual $p$, permutations $\sigma\_p$
of alternatives such that for every alternative $k$,
$a_{ik}^p\ge a_{j\,\sigma\hmii\_p\hmii(k)}^{\sigma(p)}$. This
definition is not perfectly constructive, and the authors give an
algorithm for obtaining this partial ordering, characterize it (see
Table~\ref{TCha}) and thoroughly study its properties. The partial
ordering of permuted dominance is not generated by any scoring
procedure but it is the intersection of all point orderings. We present
it here mainly because its definition is related to our
self-consistency axiom (Section~\ref{Self}).

{\em Generalized row sums\/} make up a parametric family of scores that
coincide with the extended Borda scores on complete preference
structures and generally satisfy specific systems of linear equations.

Finally, {\em Kemeny median\/} consists of all social orderings that are
nearest to the individual preferences in the metrics `sum of absolute
differences at all pairs of alternatives'. This ranking procedure is
not based on scores, and we included it because its characterizations
are closely related to those of scoring methods. As Young
\cite{Young86} revealed, this median approach had been initially
proposed in a vague form by Condorcet.
\bigskip

Thus, the last column of Table~\ref{TCha} induces the following
classification of the axiomatic characterizations. The scoring, ranking
and choice procedures in Table~\ref{TCha} are based on:
\begin{itemize}
\item Various manifestations of the extended Borda scores
\cite{Mork71,Gard73,FinFin74,Young74,HanSah,FarNit79,Barb,Rubi,NitRub,%
Henr,Young86,SaariSCW,BarGol,Debo,Bouy92a,Bouy92b,Marc96a,Marc96b}.
\item Generalized row sums \cite{Cheb89a,Cheb94}.
\item Down-sided, up-sided and factored Borda scores for profiles of
weak orders and valued relations:
   \begin{itemize}
   \item Down-sided Borda scores \cite{Dugg}.
   \item Meet of down-sided and up-sided Borda orderings \cite{BouPer}.
   \item Factored Borda scores \cite{Mork71,Mork82}; scores of approval
         voting \cite{Fish78,Fish79,BaiXu91}.
   \end{itemize}
\item Point scores and lobby size scores:
   \begin{itemize}
   \item Point scores \cite{Smith73,FinFin74,Young75,Barb}; plurality
         scores \cite{Rich78,Mork82}.
   \item Partial ordering of permuted dominance \cite{FinFin74}.
   \item Composite point scores \cite{Smith73,Young75}.
   \item Transfinite point scores \cite{FinFin74}.
   \item Lobby size scores \cite{Barb}.
   \item Convex combinations of point scores and lobby size scores
         \cite{Barb}.
   \end{itemize}
\item Sum of arbitrary scores (utilities) that represent individual
      orders \cite{Flem52,Fish69,Fish73}.
\item Aggregation functions for individual utilities $v\_1(x)\cdc
      v\_m(x)$:
   \begin{itemize}
   \item Utilitarian function $\sum_{p=1}^m\vpx$ (counterpart of the
         extended Borda scores)
         \cite{GooMar,Miln54,dAspGev77,Mask78,DesGev78,Moul88}.
   \item Weighted utilitarian functions
         $\sum_{p=1}^m\a\_p\vpx+\beta$
         \cite{BlaGir54,Hars55,Keen76,Robe80,Myer81} ($\beta\not\equiv0$
         in~\cite{Hars55}).
   \item Generalized utilitarian functions
         ${\rm sign}\beta\sum_{p=1}^m\a\_p(\vpx)^\beta$ and\newline
         $\sum_{p=1}^m\a\_p\log(\vpx)$ \cite{Robe80}.
   \item Anonymous separable functions $\sum_{p=1}^mw(\vpx)$ with
         continuous and increasing $w$ (analog of point scores) and
         continuous separable functions $\sum_{p=1}^mw\_p(\vpx)$
         \cite{Mask78}.
   \item Homothetic and increasing functions $W(v\_1(x)\cdcc v\_m(x))$
         and special cases with homogeneous of degree 1 functions
         \cite{Robe80}, see also \cite{Yano88,Khme95}.
   \item Increasing and/or continuous or arbitrary functions
         $W(v\_1(x)\cdcc v\_m(x))$ \cite{GooMar,Mask78,Robe80}.
   \end{itemize}
\item Sum of arbitrary scores determined by individuals' ballots
      \cite{Myer}.
\item Kemeny median \cite{YouLev,Young86}.
\end{itemize}

Let us mention two other related families of procedures based on
scores, first, the {\em generalized positional voting methods\/} by
Saari \cite{SaariAOR,SaariGV}. The choice is determined by point scores
but in a more complicated (and still quite natural) way involving
multiple comparisons of score differences and score sums with
thresholds.  Ordinary choice procedures based on scores are included.
Saari characterizes this family, and the key axiom is {\em weak
reinforcement}: $C(\A)=C(\B)\Rightarrow C(\A+\B)=C(\A)=C(\B)$, where
$C(\A)$ is the choice set. Merlin \cite{Merl95} demonstrates that a
version of this axiom (and of a medium condition, {\em inclusive
reinforcement}: $C(\A)\subseteq C(\B)\Rightarrow C(\A+\B)=C(\A)$),
where $C(\A)$ is the set of chosen linear orders, holds for a family of
{\em runoff ranking procedures based on scores}. Here at each stage,
some kind of point scores (for example, Borda scores or plurality
scores) is counted for a restricted set of alternatives, i.e., for
truncated preferences. The restriction of the set of alternatives is
done by eliminating from consideration low scoring (or high scoring)
alternatives which are thereby ordered. Characterizations of {\em
runoff choice procedures\/} by Hare, Coombs and Nanson are given in
\cite{Merl94}.
\label{SSResu}

\subsection{Types of preference aggregation procedures}

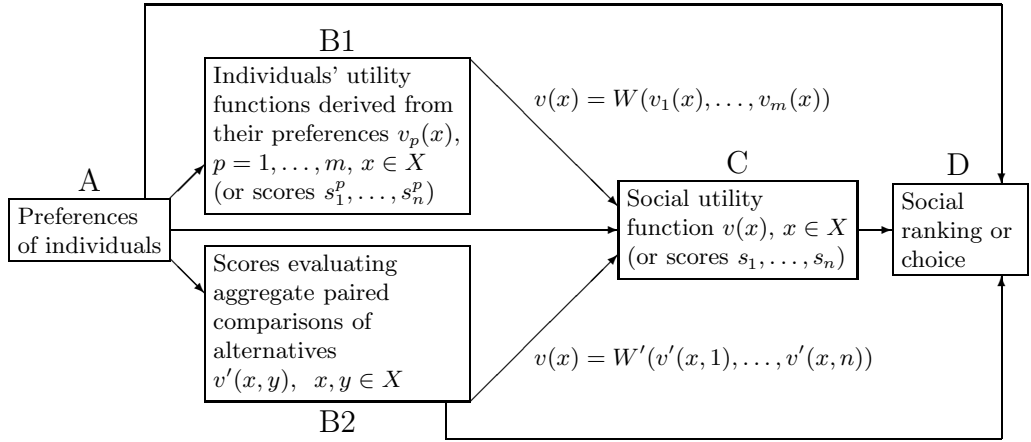
\begin{figure}[htb]
\begin{center}
\setlength{\unitlength}{5mm}
\begin{picture}(27,11.0)
\put(0,1.0){\begin{picture}(27,8.8)
\put(0, 3.7){\fbox{\parbox{19.1mm}{\raggedright Preferences of
                               individuals}}}
\put(5.2, 6.2){\fbox{\parbox{33mm}{\raggedright
                                Individuals' utility functions
                                derived from their preferences
                                $v\_p(x),$ $p=\om,$ $x\in X$ (or scores
                                $s_1^p\cdc s_n^p$)}}}
\put(5.2, 1.2){\fbox{\parbox{33mm}{\raggedright
                                Scores evaluating aggregate paired
                                comparisons of alternatives $v'(x,y),
                                \;\;x,y\in X$}}}
\put(16.2,3.7){\fbox{\parbox{29.5mm}{\raggedright
                                Social utility function $v(x)$,
                                $x\in X$ (or scores $s\_1\cdc
                                s\_n$)}}}
\put(23.5, 3.71){\fbox{\parbox{15.7mm}{\raggedright
                                Social ranking or choice}}}
\put(1.84, 4.90){\large A}            
\put(8.24, 8.62){\large B1}           
\put(8.24,-1.52){\large B2}           
\put(19.08,5.37){\large C}            
\put(24.94,5.32){\large D}            
\put(4.30,3.85){\vector(1, 0){11.92}} 
\put(4.30, 4.7){\vector(1, 1){.91}}   
\put(4.30,3.09){\vector(1,-1){.91}}   
\put(3.62,4.69){\line(0, 1){5.19}}    
\put(3.62,9.86){\line(1, 0){22.78}}   
\put(26.4,9.86){\vector(0,-1){4.77}}  
\put(11.62,-0.71){\line(0,-1){1.00}}  
\put(11.62,-1.70){\line(1, 0){14.78}} 
\put(26.4,-1.70){\vector(0,1){4.35}}  
\put(12.29,8.40){\vector(1,-1){3.91}} 
\put(12.29,-0.70){\vector(1,1){3.91}} 
\put(22.55,3.85){\vector(1,0){.94}}   
\put(13.98,7.17){$v(x)=W(v\_1(x)\cdc v\_m(x))$}
\put(13.98,.28){$v(x)=W'(v'(x,1)\cdc v'(x,n))$}
\end{picture}}
\end{picture}
\end{center}
\caption{An overall scheme of preference aggregation}
\label{Block}
\end{figure}

An overall scheme of preference aggregation is depicted in
Fig.~\ref{Block}. Block `C' is the main common feature of the scoring
methods presented in Table~\ref{TCha}. Indeed, the very term `scoring'
means `calculating some numerical indices of performance'. The methods
vary in their starting point (A or B1 or B2), destination (C or D) and
route between them.

\footnotesize
\normalsize

\bigskip{\bf A\mto B1\mto C\mto D}.
This is a prevalent route \cite{Flem52,Fish69,Mork71,Smith73,Fish73,%
FinFin74,Young75,Fish78,Rich78,Mork82,BaiXu91,Myer}. The most general
result is that given by Myerson \cite{Myer} where individual ballots
are members of an arbitrary nonempty and finite set. They even need not
be structures on the set of conceivable alternatives. For example, they
may be some colors (red, rose, brown, green, etc.) representing various
political orientations. Each ballot induces some scores of the
available alternatives (these scores can be assigned by the planner);
then the scores are summed up over individuals to give the ultimate
scores of alternatives. Anonymity and some kind of neutrality are
assumed.

Other highly general results are provided by Fleming \cite{Flem52} and
Fishburn \cite{Fish69,Fish73}. Here, individual preferences are
represented by binary relations, and the social preference structure is
determined by the sums of {\em arbitrary\/} scores (utilities) that
monotonically represent the individual relations. From the beginning,
neither neutrality nor anonymity is imposed, but then Fishburn
\cite{Fish73} studies the impact of these conditions in the framework
of social choice functions.

Point scoring methods are the neutral and anonymous variant of such
procedures. They have transfinite and composite extensions
\cite{Smith73,FinFin74,Young75}.

In \cite{Mork71,Fish78,Rich78,Mork82,BaiXu91} specific scoring methods
including plurality choice, approval voting and factored Borda method
are characterized.

\bigskip{\bf A\mto B\mto C\mto D}.
Definition (\ref{EBS}) of extended Borda scores involves two sums
which can be written in either order. Consequently, the corresponding
scoring methods (as well as those based on down-sided or up-sided
Borda scores) may take either {A\mto B1\mto C\mto D} or {A\mto B2\mto
C\mto D} route. {A\mto B\mto C\mto D} is their common designation. Such
scoring methods are characterized in
\cite{Mork71,Gard73,FinFin74,Young74,HanSah,FarNit79,Rubi,NitRub,Henr,%
Young86,SaariSCW,BarGol,Debo,Bouy92a,Bouy92b,BouPer,Dugg,Marc96a}.
In some of these papers, the input preference structure is single
binary relation (or weighted relation). It can be thought of as
representing the preferences of a single individual (the corresponding
scheme {(A$=$B1)\mto (B2$=$C)\mto D} results when `B' is divided into
`B2' and `B1' in {A\mto B\mto C\mto D}) or as an aggregate preference
relation, e.g., majority relation. The latter route is

\bigskip{\bf B2\mto C\mto D},
and the papers are \cite{Rubi,Henr,BarGol,Bouy92a,Bouy92b,BouPer}.

\bigskip{\bf B1\mto C\mto D}.
This route (A$=$B1 is usually implied) is typical of welfare economics
and game theory. The procedures characterized in
\cite{GooMar,BlaGir54,Miln54,Hars55,dAspGev77,Mask78,DesGev78,Robe80,%
Myer81,Moul88,Yano88,Khme95}
are listed in subsection~\ref{SSResu}. For other related results and
a more general context we refer to \cite{Sen77b,Dasp85,Fish87,Moul88}.

\bigskip{\bf A\mto B1\mto C, A\mto B2\mto C and A\mto B1,B2\mto C}
are represented by the characterizations of point scores, lobby size
scores and their convex combinations in \cite{Barb}. The generalized
row sum method \cite{Cheb89a,Cheb94} allows an {A\mto B1,B2\mto C}
representation (as well as many other indirect scoring procedures; some
of them---in versions suitable for complete preferences---are presented
in Table~\ref{TInd}).

\bigskip{\bf A\mto B\mto C}
corresponds to the extended Borda scores \cite{Barb,Marc96b}.

\bigskip{\bf A\mto C}.
This is for indirect scoring methods whose ultimate scores cannot be
represented through B1 and B2 structures.

\bigskip{\bf B1\mto C}
applies to \cite{Keen76}.

\bigskip{\bf A\mto B2\mto D}.
Kemeny median characterized in \cite{YouLev,Young86} does not need
calculating any social scores (`C'). This method is presented in
Table~\ref{TCha} since its characterization has much common with those
of scoring methods.

\bigskip{\bf B2\mto D}.
We refer to the recent monograph \cite{Lasl} for a comprehensive study
of such methods. The approach `based on social binary comparisons of
alternatives' is compared with the `positional' approach in
\cite{Fish71,Sen77a,Gard73,FinFin74,SaariGV}.
\bigskip

The last issue we touch on in this section is: How do the presented
characterizations of scoring methods that operate on individual orders
and individual utilities help characterize indirect scoring procedures
intended for arbitrary paired comparisons?

First, as \cite{NitRub,Young86,Debo} demonstrate, the classical
Young's characterization of the Borda method can be adapted for
arbitrary paired comparisons. The way of exploring further variations
of the axioms involved seems very promising.

Further, note that the numerical paired comparisons have the same
relation to individual utilities as two-valued paired comparisons do to
individual orders. This clears the way of interpreting and adjusting
the results of utility aggregation for the paired comparison context.

Another approach is suggested by Myerson's \cite{Myer} result. In
case we accept the conditions of Myerson's theorem, we have to admit
that the aggregation of individual paired comparisons must be based on
summing up some utilities derived from the ballots. Then the only
problem remains to introduce proper individual scores that represent
paired comparisons. The corresponding characterization problems seem to
be less complicated than the initial ones.

It is worth mentioning, however, that no scoring procedure that
operates on incomplete paired comparisons and satisfies
self-consistency, neutrality and anonymity can be represented as
{A\mto B1\mto C} \cite{CheSha97B}.
\bigskip

The procedures we characterize in the reminder of the paper operate on
profiles of valued and ordinary binary relations and result in scores
(A\mto C type). We use a unique axiom called {\em self-consistency}. It
is introduced in the following section and belongs to the `positive
relation' group (its version applicable to ranking procedures was
explored in \cite{CheSha97}). The scoring procedures that satisfy
self-consistency are turned out to have a {\em monotone implicit
form\/} (Sections~\ref{MIF}, \ref{Theo}).

\section{Self-consistency}
\label{Self}

If the scores generated by a scoring procedure are intended to serve as
numerical estimates of the alternatives, they should be comparable.
This means that whenever we consider $i$ to `perform better' than $j$,
the score of $i$ should be greater, whether they are taken from the
same or from different profiles.

Let a scoring procedure $\w$ be fixed, so the score vector $(s\_1\cdc
s\_n)$ is determined for every $\A\in\AA$.

The idea of the following axiom is simple. If $i$ got better comparison
results than $j$ and the opponents of $i$ are assigned respectively
higher score than those of $j$, then the score of $i$ should be greater
than the score of $j$. To formalize this requirement, we introduce some
special notation. Recall that $M=\{\om\}$ and
$X\_i=\{\on\}\setminus\{i\}$, $i=\on$ (Section~\ref{Nota}).

\begin{definition}[(Performance multiset)]
The multiset\footnote{%
For multisets, as distinct from sets,
multiple occurrence of elements is allowed. $U\_i$ is a multiset
since the numerical pairs $(a_{ik}^p, s\_k)$ may coincide for
different $k$ and $p$
}
$U\_i=\{(a_{ik}^p, s\_k)\mid k\in X\_i,\:p\in M\}$ of pairs $(a_{ik}^p,
s\_k)$ corresponding to all comparison outcomes of $i$ in $\A$ will be
referred to as the {\em performance multiset} of $i$ in $\A$.
\end{definition}

Let $\A$ and $\A'$ be two admissible profiles of individual
preferences.  Suppose $U\_i$ and $U'_j$ are the performance multisets
of $i$ in $\A$ and $j$ in $\A'$, respectively.

\begin{definition}[(Majorization)]
Alternative $i$ in $\A$ {\em majorizes\/} $j$ in $\A'$ if there exists
a one-to-one mapping $\mu$ from $U\_i$ onto $U'_j$ such that
$\mu\Big((a_{ik}^p,s\_k)\Big)=(a'^q_{j\l},s'_{\l})$ implies
$a_{ik}^p\ge a'^q_{j\l}$ and $s\_k\ge s'_{\l}.$ Furthermore, $i$ in
$\A$ {\em strictly majorizes\/} $j$ in $\A'$ if, in addition, at least
one of the above inequalities is strict for at least one comparison
outcome $a_{ik}^p$.
\end{definition}

\begin{axiom}[(Self-consistency)]
1. If $i$ in $\A$ majorizes $j$ in $\A'$ then $\si\ge s'_j$.\\
2. If $i$ in $\A$ strictly majorizes $j$ in $\A'$ then
$\si>s'_j$.
\end{axiom}

It can be said that a scoring procedure is self-consistent when it
{\it preserves majorization}. In other words, this axiom is a kind of
Pareto condition with a self-consistent version of superiority which
allows permutations and inter-profile juxtapositions. Among the axioms
of the previous section, self-consistency is related with the Pareto
principle for equicardinal subsets \cite{Fish69,Fish73}, permuted
dominance \cite{Fish71b} and especially with its version exploited in
\cite{FinFin74}. As distinct from them, self-consistency applies to
A\mto C aggregation procedures, which enables it to capture
interprofile comparisons. Besides, self-consistency is weaker in that
it only recognizes the superiority in comparison outcomes confirmed by
exceeding scores of the `opponents'.

\begin{figure}[htb]
\begin{center}
\setlength{\unitlength}{5mm}
\begin{picture}(21.9,11.6)
\put( 5, 0){\line(1,0){12}}
\put( 5,11){\line(1,0){12}}
\put( 5, 0){\line(0,1){11}}
\put(17, 0){\line(0,1){11}}
\put( 5, 0){\begin{picture}(8,11)
    \put( -3, .5){\begin{picture}(12,11)
             \put(4.0,6.0){$i$}
             \put(8.0,3.0){$e$}
             \put(8.0,6.0){$c$}
             \put(8.0,9.0){$a$}
             \put(7.7,8.9){\vector(-4,-3){3.1}} 
             \put(4.5,6.1){\vector(1  ,0){3.0}} 
             \put(4.5,5.7){\vector(4, -3){3.0}} 
             \put(6.0,.2){$\A$}                 
    \end{picture}}
    \put( -1.2, .1){\begin{picture}(8,11)
             \put(8.0,2.0){$f$}
             \put(8.0,5.0){$d$}
             \put(8.0,8.0){$b$}
             \put(12.0,5.0){$j$}
             \put(8.5,7.9){\vector(4 ,-3){3.2}}  
             \put(11.7,5.1){\vector(-1 ,0){3.0}} 
             \put(8.6,2.2){\vector (4,  3){3.1}} 
             \put(10.0,.6){$\A'$}                
    \end{picture}}
    \multiput(6.1,0.0)(0,1.04){11}{\line(0,1){.6}}
    \end{picture}}
\end{picture}
\end{center}
\caption{An illustration to self-consistency}
\label{SCIll}
\end{figure}
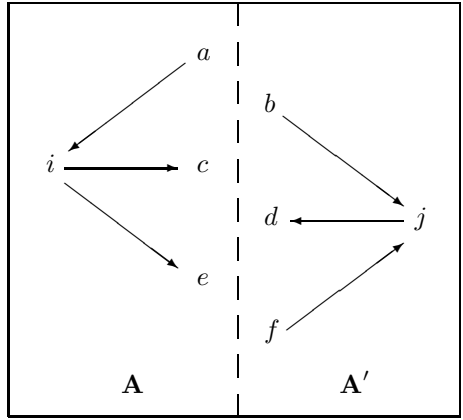

Self-consistency is illustrated in Fig.~\ref{SCIll} where only
comparison outcomes of $i$ and $j$ are shown.
$U\_i=\{(0,s\_a),(1,s\_c),(1,s\_e)\}$ and
$U'_j=\{(0,s'_b),(1,s'_d),(0,s'_f)\}$.
If $s\_a\ge s'_b$, $s\_c\ge s'_d$, and $s\_e\ge s'_f$, then $i$ in
$\A$ strictly majorizes $j$ in $\A'$ with the following $\mu$:
$\mu\Big((0,s\_a)\Big)=(0,s'_b)$,
$\mu\Big((1,s\_c)\Big)=(1,s'_d)$,
$\mu\Big((0,s\_a)\Big)=(0,s'_f)$.
Then self-consistency requires $s\_i>s'_j$.

\section{Monotone implicit form}
\label{MIF}

Now consider scoring procedures of different nature. They came from
such disciplines as management science, psychometrics, applied
statistics, processing of sport tournaments, graph theory, etc., and
are based on the resolution of systems of algebraic equations. The
number of equations is the number of alternatives and the form of them
is presented for several procedures in Table~\ref{TInd}. Five of them
were rediscovered for several times and with different motivations. In
the first column we give only the earliest references we know, other
ones can be found in \cite{CheSha96,CheSha97B}. In Table~\ref{TInd}
these procedures are adjusted to the type of data considered here. The
domain of them not necessarily contains all possible profiles, but this
subject is out of our scope now. The method by Smith and Gulliksen was
designed for incomplete preference data.

\begin{Table*}{3}
{Some sensitive scoring procedures}
\label{TInd}
Paper&
$i$th equation of $n$ equations&
Requirements or corollaries\\
\hline
\noalign{\vspace{2pt}}
Zermelo \cite{Zerm}&
$\suml_{j\ne i}\suml_{p=1}^m
\left(a_{ij}^p-\frac{\textstyle\si}{\textstyle\si+\sj}\right)=0$&
$\si>0,\:\suml_{i=1}^n \si=1$\\
\noalign{\vspace{2pt}}
Katz \cite{Katz}&
$\suml_{j\ne i}\suml_{p=1}^m
a_{ij}^p\left(\ve\sj+1-\frac{\textstyle\si}{\textstyle
m(n-1)}\right)=0$&
$\ve>0$, $\si>0$\\
\noalign{\vspace{3pt}}
\begin{tabular}{l}
Smith \cite{Smith56},\\
Gulliksen \cite{Gull}\\
\end{tabular}&
$\suml_{j\ne i}\suml_{p=1}^m
\big(mn(a_{ij}^p-a_{ji}^p)+\sj-\si\big)=0$&
$\suml_{i=1}^n \si=0$\\
\noalign{\vspace{3pt}}
Daniels \cite{Dani}&
$\suml_{j\ne i}\suml_{p=1}^m
(a_{ij}^p\sj-a_{ji}^p\si)=0$&
$\si>0$\\
\noalign{\vspace{2pt}}
Daniels \cite{Dani}&
$\suml_{j\ne i}\suml_{p=1}^m
\left(a_{ij}^p\frac{\textstyle\sj}{\textstyle\si}-
a_{ji}^p\frac{\textstyle\si}{\textstyle\sj}\right)=0$&
$\si>0$\\
\noalign{\vspace{2pt}}
Cowden \cite{Cowd}&
$\suml_{j\ne i}\suml_{p=1}^m
(a_{ij}^p\sj(1-\si)-a_{ji}^p\si(1-\sj))=0$&
$\si>0$\\
\noalign{\vspace{2pt}}
\end{Table*}

The common idea of these procedures is to process comparison outcomes
taking into account the strength of competing alternatives which
(strength) is estimated through their `performance' in the same
profile. In this sense, these methods are self-consistent. A relation to
the axiom of the same name is conveyed by the following theorem.

It is easily seen that each equation in Table~\ref{TInd} has the form
\begin{equation}
\suml_{j\ne i}\suml_{p=1}^m h(a_{ij}^p,\sj,\si)=0,\qquad i=\on,
\label{impsi}
\end{equation}
where $h(\cdot,\cdot,\cdot)$ strictly increases in $a_{ij}^p$ (recall
that $a_{ij}^p+a_{ji}^p=1$, $j\ne i$, $p=\om$, as assumed in
Section~\ref{Nota}) and in $\sj$ and strictly decreases in $\si$.
These properties of $h(\cdot,\cdot,\cdot)$ can be intuitively motivated
as follows: The greater comparison outcomes $a_{ij}^p$ and `strengths
of opponents' $\sj$ alternative $i$ has, the greater should be its own
`strength' $\si$ (to provide the zero sum on the left-hand side).

A form close to (\ref{impsi}) was used as an axiom to derive the
generalized row sum method in \cite{Cheb89a,Cheb94}. It can be
noted however that this form is somewhat too special because of the
double sum on the left-hand side. Let us replace this sum with an
arbitrary strictly increasing function.

Let $\T$ be the set of all multisets that consist of $m(n-1)$ real
triples.

\begin{definition}[(Monotone implicit form of scoring procedure)]
A scoring procedure $\w:\:\AA\to \R^n$ has a {\em monotone
implicit form\/} if there exists a function $g:\T\to\R$ such that\\
{\rm (i)} for every profile $\A\in\AA$, the scores satisfy the
system of equations
\begin{equation}
g\big(\{(a_{ij}^p,\sj,\si)\mid j\in X\_i,\:p\in M\}\big)=0,
\qquad i=\on;
\label{MIFe}
\end{equation}
{\rm (ii)} $g$ strictly increases in every $a_{ij}^p$ and $\sj$
and strictly decreases in $\si$.
\end{definition}

This form is a direct generalization of (\ref{impsi}).

\section{Self-consistency amounts to the existence of a monotone
implicit form}
\label{Theo}

\begin{theorem}
A scoring procedure $\w$ is self-consistent if and only if it has a
monotone implicit form.
\label{Th}
\end{theorem}

The `if' part of the theorem is an easy consequence of the
corresponding definitions (see the proof). The converse statement is
not so trivial. Below it is reduced to the fact that every bounded
function defined on any `Paretian subset' of $\R^k$ has a strictly
monotonic extension to the $\R^k$ (Lemma~\ref{Lem}). In fact, we need a
special case of Lemma~\ref{Lem}, but its general formulation has
essentially the same proof and is worth mentioning by itself.

Prior to proving Theorem~\ref{Th}, note that the extended Borda scores
(\ref{EBS}) trivially satisfy self-consistency (the proof is left to
the reader). By Theorem~\ref{Th}, this procedure must have a monotone
implicit form. Indeed, as shown in \cite{Cheb89a,Cheb94}, a
one-parametric family of such forms is provided by the generalized row
sum method:
\begin{equation}
\suml_{j\ne i}\suml_{p=1}^m(\gamma(a_{ij}^p-a_{ji}^p)-
(\si-\sj))-\si\ve^{-1}=0,\qquad i=\on,
\label{GRSf}
\end{equation}
where $\ve>0$ is a positive parameter and $\gamma=\ve^{-1}+mn$. The
least squares procedure by Smith \cite{Smith56} and Gulliksen
\cite{Gull} presented in Table~\ref{TInd} provides another form which
can be obtained from (\ref{GRSf}) as $\ve\to\infty$.

\begin{proof}
Let $\w$ be a scoring procedure such that there exists $g:
\T\to\R$ having the properties (i) and (ii) of monotone
implicit form.  Prove that $\w$ is self-consistent.

Suppose that $i$ in $\A$ majorizes $j$ in $\A'$. Then there
exists a one-to-one mapping $\mu$ from $U\_i$ onto $U'_j$ such that
$\mu\Big((a_{ik}^p,s\_k)\Big)=(a'^q_{j\l},s'_{\l})$ implies
$a_{ik}^p\ge a'^q_{j\l}$ and $s\_k\ge s'_{\l}.$ Assume that
$\si<s'_j$. Then by (ii) the left-hand side of the $i$th
equation of (\ref{MIFe}) written for $\A$ is greater than the
left-hand side of the $j$th equation of (\ref{MIFe}) written for
$\A'$ and they cannot be both equal to zero. Hence $\si\ge s'_j$
and item 1 of self-consistency is satisfied. Item 2 is proved
similarly.

Now suppose that $\w$ is self-consistent. Prove that there exists
a function $g: \T\to\R$ satisfying (i) and (ii).

Let us say that a set $P\subset \R^k$, $k\in\N=\{1,2,\ldots\}$ is
a {\it Paretian subset of\/} $\R^k$ if for any $(\z,\z')\in P^2$,
either $\z=\z'$ or there exists $i\in \{\ok\}$ such that
$z_i>z'_i$.

For every profile $\A\in\AA$ and every alternative $i$, the multiset
$\{(a_{ij}^p,\sj,-\si)\mid j\in X\_i,\:p\in M\}$, where
$a_{ij}^p$ are the comparison outcomes of $i$ in $\A$ and
$\si,\:\sj$ are the scores assigned by $\w$, will be called the
{\it multiset of comparison triples of $i$ in $\A$}.

Let $P\_{\w}\subset \R^t$ where $t=3m(n-1)$ be the set of all
vectors $\z=(z\_1,z\_2\cdc z\_t)$ such that the multiset $\{(z\_1,z\_2,
z\_3),(z\_4,z\_5, z\_6)\cdc(z\_{t-2},z\_{t-1}, z\_t)\}$ is the
multiset of comparison triples of some alternative in some profile
$\A\in\AA$. It follows from self-consistency that $P\_{\w}$ is a
Paretian subset.
\bigskip

\begin{lemma}
\label{Lem}
Suppose $P$ is a Paretian subset of \/
$\R^k$, $k\in\N$. For every bounded function $f\_P(\x):\;
P\to\R$, there exists a function $f(\x):\;
\R^k\to\R$ such that
\par{\rm($\ast$)} The restriction of $f(\x)$ to $P$ coincides with
$f\_P(\x)$ and
\par{\rm($\ast\ast$)} $f(\x)$ is strictly increasing in every component of
$\x${\rm:} $x\_1\cdc x\_k$.
\end{lemma}

\begin{proof}
1. Prove Lemma~\ref{Lem} for the open hypercube\footnote{$i$, $j$ and
$k$ in the proof of Lemma~\ref{Lem} have nothing in common with the
same variables that denote alternatives} $L=\Big\{\x=(x\_1\cdc
x\_k)\;\Big|\; |x\_i|<1,\; i=\ok \Big\}$ substituted for the Euclidean
space $\R^k$. After that the general case of $\R^k$ will be reduced to
this one.

Introduce the following notation.

For any $\x,\y\in L,\:$ $\x\ge \y$ means $\forall i\in\{\ok\}\:$
$x\_i\ge y\_i$.

For every $\x\in L,$ let
$$
D(\x)=\Big\{\y\in L\;\Big|\; \x\ge \y \Big\}\cup(\R^k\setminus L),
$$
$$
U(\x)=\Big\{\y\in L\;\Big|\; \y\ge \x \Big\}\cup(\R^k\setminus L).
$$

For any $Y\subseteq L$, let
$$
D(Y)=\cupl_{\x\in Y} D(\x),
$$
$$
U(Y)=\cupl_{\x\in Y} U(\x).
$$

Suppose $\e^i$ is the unit vector of the $i$th axis of $\R^k$;
$D=D(P),$ $U=U(P)$.  Suppose that $F\_{\min}, F\_{\max}\in\R$ are
such that for any $\x\in P$, $F\_{\min}\le f\_P(\x)\le F\_{\max}$
holds (recall that $f\_P(\x)$ in Lemma~\ref{Lem} is bounded). Define
several auxiliary functions on $L$.
$$
d^i(\x)=\inf\Big\{d\ge0\;\Big|\;\x-d \e^i\in D\Big\},\quad i=\ok,
$$
$$
u^i(\x)=\inf\Big\{u\ge0\;\Big|\;\x+u \e^i\in U\Big\},\quad i=\ok.
$$

For every $\x\in L$ define
$$
f\_1(\x)=\suml_{i=1}^k d^i(\x)-\suml_{i=1}^k u^i(\x),
$$
$$
f\_2(\x)=\cases{F\_{\min}, &$\x\in L\cap D\setminus U$  \cr
               F\_{\max}, &$\x\in L\cap U\setminus D$  \cr
               {1\over 2}(F\_{\min}+F\_{\max}),
                          &$\x\in L\setminus(D\cup U)$ \cr
               f\_P(\x),   &$\x\in L\cap D\cap U=P,$    \cr
}
$$
and finally,
$$
f(\x)=f\_1(\x)+f\_2(\x).
$$

Now we prove that $f(\x)$ possesses the desired properties ($\ast$) and
($\ast\ast$). First, note that $d^i(\x)$ and $u^i(\x)$, $i=\on$, and thus
$f(\x)$ are well-defined by the definitions of $D$ and $U$.

($\ast$) For any $\x\in P$, $f(\x)=f\_P(\x)$ since $\forall \x\in P$
$\;d^i(\x)=u^i(\x)=0$, and $f\_2(\x)=f\_P(\x)$.

($\ast\ast$) Prove that $f(\x)$ is strictly increasing on $L$ in every
$x\_i$.  It suffices to show that for any $\x\in L$, $\a>0$, and
$i\in\{\ok\}$,
$$
\x'=\x+\a \e^i\in L\mbox{\quad implies\quad} f(\x')>f(\x).
$$

Fix $i$ and prove three statements.

(A) $d^i(\x')-u^i(\x')>d^i(\x)-u^i(\x).$

It is easily seen that $d^i(\x')\ge d^i(\x)$ and $u^i(\x')\le
u^i(\x)$. Assume that $d^i(\x')=d^i(\x)$ and $u^i(\x')=u^i(\x)$. Then
$\x,\x'\in L\cap D\cap U=P$. This is impossible since $P$ is
Paretian.
\medskip

(B) For every $j\in \{\ok\}\setminus \{i\}$, $d^j(\x')\ge d^j(\x)$ and
$u^j(\x')\le u^j(\x)$.

Assume that $d^j(\x')<d^j(\x).$ Then there exists $d\_0\ge0$ such
that $\x'-\dej\in D$ and $\x-\dej\notin D$. According to the definition
of $D$, this is impossible, since if there exists $\z\in P$ such that
$\z\ge \x'-\dej$, then $\z\ge \x-\dej$ and hence $\x-\dej\in D$; on the
other hand, if $\x'-\dej\in \R^k\setminus L$, then $\x-\dej\in
\R^k\setminus L$ and $\x-\dej\in D$ as well. Hence $d^j(\x')\ge d^j(\x)$.
Similarly, $u^j(\x')\le u^j(\x)$.
\medskip

(C) $f\_2(\x')\ge f\_2(\x)$.

It is easy to verify that all possible translations of the vector
$\x+\a \e^i\in L$ from one set to another as $\a\ge0$ increases
are specified by the following diagrams:
$$
\begin{array}{llcrlr}
L\cap D\setminus U&\to &P=L\cap D\cap U
             &\to &L\cap U\setminus D &\mbox{ or}\\
L\cap D\setminus U&\to &L\setminus(D\cup U)
             &\to &L\cap U\setminus D &\mbox{ or}\\
L\cap D\setminus U&\; &\to&\; &L\cap U\setminus D.&
\end{array}
$$
At all these translations, $f\_2(\x)$ does not decrease by
definition.
\medskip

By (A), (B) and (C), $f(\x)$ is strictly increasing in $x\_i$ on $L$.
\medskip

2. Now let $P$ be any Paretian subset of $\R^k$. To extend $f\_P(\x)$ to
$\R^k$, we first contract  $\R^k$ onto $L$ by the mapping $\y=\psi(\x)$
where $y\_i={2\over \pi} \arctan x\_i,$ $i=\ok$. Then solve the problem
on $L$ and finally extend $L$ onto $\R^k$ by $\y=\psi^{-1}(\x)$, i.e.,
$y\_i=\tan({\pi\over 2} x\_i)$, $i=\ok$. It is obvious that under these
strictly increasing transformations, Paretian subsets are mapped to
Paretian subsets and strictly increasing functions to strictly
increasing functions.  This completes the proof of Lemma~\ref{Lem}.
\end{proof}

To prove Theorem~\ref{Th}, apply Lemma~\ref{Lem} to $P=P\_{\w}$.
Namely, put $\R^k=\R^t$ and $f\_{P\_{\w}}\equiv0.$ Then the conclusion
of Lemma~\ref{Lem} differs from the desired statement only by that the
required function $g$ should be defined for the multisets of real
triples and should map all {\em multisets of comparison triples\/} to
zero, whereas function $f$ provided by Lemma~\ref{Lem} is defined on
$\R^t$ and maps to zero all vectors corresponding to the multisets of
comparison triples. But note that $P\_{\w}$ is invariant (by
definition) to any permutation of triples of adjacent coordinates
$(z\_{3r-2},z\_{3r-1}, z\_{3r})$. Then, by the proof of
Lemma~\ref{Lem}, in this case $f$ is also invariant to such
permutations applied to its argument. This implies that the value of
$f$ for any vector $\x\in\R^t$ is solely determined by the {\em
multiset\/} of adjacent triples of coordinates $x\_1\cdc x\_t$. In
other words, $f$ does not depend on the order of these triples in
$(x\_1\cdc x\_t)$.  Therefore, $f$ being considered as a function of
such multisets of triples defines a function $g$ with the desired
properties. This completes the proof.

\end{proof}

\end{document}